\documentclass[a4paper, 11pt]{article}

\usepackage{amsmath,pifont,amssymb}
\usepackage{psfrag}
\usepackage{subfigure}
\usepackage{url}
\usepackage{color}
\usepackage[latin1]{inputenc}
\usepackage[T1]{fontenc}
\usepackage[english]{babel}
\usepackage{verbatim}
\usepackage{graphics,array,float,epsfig}
\usepackage{amsfonts,graphicx,color}
\usepackage{dsfont}
\usepackage{mathrsfs}
\usepackage{stmaryrd}
\usepackage{authblk}

\def\squarebox#1{\hbox to #1{\hfill\vbox to #1{\vfill}}}

\newtheorem{theorem}{Theorem}[section]
\newtheorem{lemma}{Lemma}[section]
\newtheorem{proposition}{Proposition}[section]
\newtheorem{remark}{Remark}[section]

\newcommand{\sub}{\scriptscriptstyle}
\newcommand{\reel}{\mathds{R}}

\newcommand{\dif}{{\mathrm{d}}}

\newcommand{\tp}{^{\scriptscriptstyle T}} 
\newcommand\R{\mathbb{R}}

\newcommand\N{\mathbb{N}}

\newcommand{\transp}{^{\scriptscriptstyle T}}

\newcommand{\Exp}{\mathbb{E}}

\newcommand{\D}{{\cal D}}

\newcommand{\alf}{\scriptscriptstyle{1/2}}

\newcommand{\mbs}[1]{\ensuremath{\boldsymbol{#1}}}
\newcommand{\bcdot}{\mbs{\cdot}}

\newcommand{\Df}{\mathrm{\mathbb{D}}} 
\newcommand{\sdbt}{\sigma \dif W_t} 
\renewcommand*{\div}{\nabla\bcdot} 
\newcommand{\adv}{\bcdot\nabla} 

\begin{document}
\pagestyle{headings}
\setcounter{tocdepth}{4}
\title{A consistent stochastic large-scale representation of the Navier-Stokes equations}
\author[1]{Arnaud Debussche}
\author[2]{Berenger Hug}
\author[,2]{Etienne M\'emin}
\affil[1]{Univ Rennes, CNRS, IRMAR - UMR 6625, F- 35000 Rennes, France }
\affil[2]{Inria/IRMAR Campus de Beaulieu 35042 Rennes Cedex}

\maketitle
\begin{abstract}
In this paper we analyze the theoretical properties of a stochastic representation of the incompressible Navier-Stokes equations defined in the framework of the modeling under location uncertainty (LU). This setup built from a stochastic version of the Reynolds transport theorem incorporates a so-called transport noise and involves several specific additional features  such as a large scale diffusion term, akin to classical subgrid models, and a modified advection term arising from the spatial inhomogeneity of the small-scale velocity components. This formalism has been numerically evaluated in a series of studies with a particular interest on geophysical flows approximations and data assimilation. In this work we focus more specifically on its theoretical analysis.  We demonstrate, through classical arguments, the existence of martingale solutions for the stochastic Navier-Stokes equations in LU form. We show they are pathwise and unique for 2D flows. We then prove that if the noise intensity goes to zero, these solutions converge, up to a subsequence in dimension $3$, to a solution of the deterministic Navier-Stokes equation.  similarly to the grid convergence property of well established large-eddies simulation strategies, this result allows us to give some guarantee on the interpretation of the LU Navier-Stokes equations as a consistent large-scale model of the deterministic Navier-Stokes equation. 
  
\end{abstract}

\section{Introduction}
For several years there has been a burst of activity to devise stochastic representations of fluid flow dynamics. These models are strongly  motivated in particular  by climate and weather forecasting issues \cite{Berner-17,Franzke2015,Gottwald2017, Majda99} and the need to provide accurate likely scenarios with proper uncertainty quantification as well as the necessity to build efficient methods for the coupling of data of ever increasing resolution with large-scale geophysical dynamical models. In the past decades, several schemes have been devised in that prospect. Multiplicative random forcing and randomization of parameters based on early turbulence studies on energy backscattering across scales \cite{Leith90,Mason92} have been  proposed \cite{Buizza99,Shutts05} for weather forecasting. These  schemes  have not been designed within clear mathematical and physical derivation settings and consequently lack generality. They theoretically face uncontrolled  variance increase and are hence mathematically not well posed. Due to this, a proper tuning of an eddy viscosity  term must be realized to counterbalance the energy brought by the noise. The precise form that should take this additional eddy viscosity term remains unknown and often relies on debatable physical hypothesis such as the Boussinesq's assumption \cite{Schmitt:2007}.  More problematically, random forcing defined outside any conservation principles may lead, even for low noise, to strong discrepancies with the fine resolution dynamics that is meant to be emulated  \cite{Chapron-18}. This issue highlights indeed that two somewhat opposite goals are  sought by stochastic parameterization. On the one hand, one wishes to augment the variability of coarse resolution systems, which is particularly desirable for data assimilation applications and uncertainty quantification, and on the other hand, to represent, as accurately as possible,  coarse resolution of the geophysical dynamics.  This latter goal is a large eddies simulation (LES) issue, while the first one seeks to avoid the traditional over-diffusive approaches associated to  coarse scale simulations embedding classical subgrid models. A minimal mathematical requirement for satisfactory LES  is that a weak solution of the LES scheme converges in some sense toward a weak solution of the fine-scale deterministic Navier-Stokes equations in 3D and toward the unique solution for the 2D Navier-Stokes equations.  The convergence of some classical LES models toward the true fine scale dynamics is well known in the deterministic case \cite{Berselli-et-al-2010, Guermond:2004}. However, the question of convergence of stochastic parametrization toward solutions of the deterministic equations at the limit of vanishing noise is not always clear. 

Recently two different general modelling schemes have been proposed to derive systematically stochastic versions of fluid dynamical models. The first one is based on a Hamiltonian principle and enables to exhibit a circulation preserving representation of Euler equations \cite{Holm-15}. This scheme has been analyzed in vorticity form \cite{Crisan-Flandoli-Holm-18} and has been assessed in several geophysical configuration approximations \cite{Cotter-MMS-19}. The second scheme, based also on stochastic transport, relies more directly on Reynolds transport theorem \cite{Memin14}. This scheme, constructed from the classical physical conservation laws preserves energy. It has been successfully applied  to several geophysical models \cite{Bauer-et-al-JPO-20, Resseguier-GAFD-I-17, Resseguier-GAFD-II-17, Resseguier-GAFD-III-17}, to the design of stochastic reduced order models \cite{ Resseguier-JFM-17, Resseguier-SIAM-UQ-21, Tissot-JFM-21} and to devise large eddies simulation models \cite{Chandramouli-CF-18, Chandramouli-JCP-20, Kadri-CF-17}. Beyond its  construction the physical relevance of the LU scheme has been numerically assessed on several prototypical flow models \cite{Bauer-et-al-JPO-20,Brecht-et-al-2021,Chapron-18}. The objective is here to check its mathematical consistency in terms of its well posedness and in terms of its ability to provide in some sense a converging stochastic representation of the Navier Stokes equations.  

The two stochastic schemes based on transport noise, proposed in \cite{Holm-15,Memin14}, are meant to represent the flow at large scales and rely on a  decomposition of the Lagrangian velocity in terms of a smooth-in-time component and a strongly oscillating component. Nevertheless the question whether or not those schemes converge towards the deterministic Navier-Stokes solutions has not been yet answered. In this work we specifically focus on the second family of methods, referred to as, modelling under location uncertainty (LU) and provide an affirmative answer to this question. We show that LU Navier-Stokes equations admit martingale solutions (also called weak probabilistic solutions) in 3D and a unique strong solution - in the probabilistic sense - in 2D. Moreover, these solutions converge toward weak solutions of the Navier-Stokes equation in 3D and toward the unique solution in 2D. 
As such these results enable to consider the LU representation as a valid large-scale stochastic representation of flow dynamics that is more amenable to ensemble forecasting and data assimilation than deterministic model due to an improved variability. The Navier-Stokes equation with transport noise has been the object of many articles, starting with \cite{BrCaFl92,Mikulevicius04}. However,
none have considered the equations studied here and the limit of the noise going to zero has not been investigated.

The paper is structured as follows. We first briefly introduce the LU modelling and the form of the corresponding stochastic Navier-Stokes equations in this setting. The main convergence results will then be enunciated. Their rigorous proofs are then fully described in the following section.  

\section{Modelling under location uncertainty}
The LU formulation relies mainly on  the following  time-scale separation assumption of the flow:
\begin{equation}\label{eq:dX}
\dif X_t = u (X_t, t)\, \dif t + \sigma (X_t, t)\, \dif W_t,
\end{equation}
where $X$ is the Lagrangian displacement defined within the bounded domain $\mathcal{S} \subset \reel^d\ (d = 2\ \text{or}\ 3)$ with smooth boundary, and  $u$ denotes the large-scale velocity that is both spatially and temporally correlated, while $\sigma\dif W$ is a highly oscillating  unresolved component (also called noise term) that is only correlated in space. 

More precisely, we consider a cylindrical Wiener process $W$ on $L^{2}(\mathcal{S} , \R^{d})$, the space of square integrable functions on $\mathcal S$ with values in $\R^d$,
$$
W=\sum_{i\in \N} \hat \beta^i e_i,
$$
where $(e_i)_{i\in \N}$ is a Hilbertian orthonormal basis of $L^{2}(\mathcal{S} , \R^{d})$ and $(\hat \beta_i)_{i\in\N}$ is a sequence of independent standard brownian motions on a stochastic basis $(\Omega , \mathcal{F},  (\mathcal{F}_t)_{t\in[0,T]}, \mathbb{P} )$ (\cite{DaPrato}). The above does not converge in $L^{2}(\mathcal{S} , \R^{d})$ but in any larger Hilbert space $U$ such that the embedding of 
$L^{2}(\mathcal{S} , \R^{d})$ into $U$ is Hilbert-Schmidt. 

The spatial structure of the noise is specified through a deterministic time dependent  integral covariance operator 
$\sigma_t$  defined from a bounded and symmetric kernel $\widehat{\sigma}$:
\[
(\sigma_t f)(x)\;  := \; \int_{\mathcal{S}} \widehat{\sigma}(x,y,t) \; f(y)\; \mathrm{d}y, \; f\in L^{2}(\mathcal{S} , \R^{d}).
\]
For each $(x,y,t)$,  $\widehat{\sigma}(x,y,t)$ is a $d\times d$ symmetric tensor. Since $\hat \sigma$ is bounded in $x,\, y$ and $t$, $\sigma_t$ maps $L^{2}(\mathcal{S} , \R^{d})$ into itself and is Hilbert-Schmidt. 
Then, the noise can be written as the Wiener process:
\[
\sigma_t W_t = \sum_{i\in \N} \hat \beta^i_t \sigma_t e_i,
\]
where the series converges in $L^{2}(\mathcal{S} , \R^{d})$ almost surely and in $L^p(\Omega)$ for all $p\in \N$ and equation \eqref{eq:dX} should be understood in the It\^o sense. To avoid, any confusion, we may further write the dependance of the Wiener process
in terms of the other variables:
\[
\sigma_t W_t (x,\omega) = \sum_{i\in \N} \hat \beta^i_t(\omega) \sigma_t e_i(x),
\]

We consider a divergence free noise:
$$
\nabla_x\cdot \hat \sigma(x,y,t)=0, \; x,y\in \mathcal  S, \; t\ge 0.
$$

Also, for each $t\in \R^+$, there exists $(\phi_n(t) )_n$ a complete orthogonal system composed of  eigenfunctions of the covariance operator at each time $t\in\mathbb{R}$. It can be shown that there exists another sequence of independent standard brownian motions, on the same  stochastic basis $(\Omega , \mathcal{F},  (\mathcal{F}_t)_{t\in[0,T]}, \mathbb{P} )$ such that we have the representation:
\[
\sigma_t \,W_t = \displaystyle \sum_{k=0}^{\infty} \phi_{k}(t) \, \beta^{k}_t.
\]

This Gaussian random field is associated to the two-times, two-points covariance tensor given by
\[
Q(x,y,t,t') = \Exp \left( \sigma_t \, \mathrm{d}W_t(x) \; \sigma_{t'} \, \mathrm{d}W_{t'}(y)\right)  = \int_{\mathcal{S}} \widehat{\sigma}(x,z,t) \, \widehat{\sigma}(z,y,t')  \mathrm{d}y \, \delta(t-t')\,,
\]
with the diagonal part (i.e one time auto-correlation), referred to,  in the following, as the variance tensor, and denoted  by
\begin{equation}
\label{def_a}
a (x,t) = \; \int_{\mathcal{S}} \widehat{\sigma}(x,y,t) \, \widehat{\sigma}(y,x,t)\; \mathrm{d}y\; = \; \displaystyle \sum_{k=0}^{\infty} \phi_{k}(x,t) \, \phi^{\scriptscriptstyle{T}}_{k}(x,t) .
\end{equation}

Some important remarks on decomposition \eqref{eq:dX} can be done at this point. It must be stressed that it is a temporal decomposition and not a spatial decomposition as classically formulated through spatial filtering and/or decimation operators in large-eddies simulation (LES) techniques \cite{Sagaut04}. However, for turbulent flows, time and spatial scales are intricately related. As a matter of fact, in the inertial range, the turn-over time ratio for two different scales $L$ and $\ell$ reads  $\tau_L/\tau_{\ell} \propto (L/\ell)^{2/3}$ and provides a direct relation between time-scale coarsening and spatial-scale dilation.  Unless specifically indicated, in the following, the resolved (unresolved) components will always be referred to as large-scales (small-scales), respectively, without differentiating between time or space scales. It can be noted that temporal filtering has already been used for the definition of oceanic models \cite{Hecht08} or large-eddies simulation approaches \cite{Meneveau00}.

In a way similar to the classical derivation of Navier-Stokes equations, the LU setting is based on a stochastic representation of the Reynolds transport theorem (SRTT) \cite{Memin14}, describing the rate of change of a random scalar $q$ within  a volume $V(t)$  transported by the stochastic flow \eqref{eq:dX}. For incompressible unresolved flows, (i.e. $\div \sigma = 0$), the SRTT reads
\begin{subequations}
\begin{align}
&\dif\, \Big( \int_{V(t)} q (x, t)\, \dif x \Big) = \int_{V(t)} \big( \Df_t q + q \div (u - u_s)\dif t \big)\, \dif x , \label{eq:SRTT} \\
&\Df_t q = \dif_t q + (u - u_s) \adv q\, \dif t + \sdbt \adv q - \frac{1}{2} \div (a \nabla q)\, \dif t, \label{eq:STO} 
\end{align}
\end{subequations}
where $\dif_t q (x, t) = q (x, t + \dif t)  - q (x, t)$ stands for the forward time-increment of $q$ at a fixed point $x$, $\Df_t$ is introduced as the stochastic transport operator in \cite{Memin14, Resseguier-GAFD-I-17} and plays the role of the material derivative. This operator is derived from the It\^{o}-Wentzell formula \cite{Kunita} to express the differentiation of a stochastic process transported by the flow \cite{Memin14}.  Recall that $u$ is the  large-scale velocity used in \eqref{eq:dX} and $a$ is defined in \eqref{def_a}. Note also that we omit to mention the dependance of $\sigma$ on time.  The drift $u_s = \frac{1}{2}\div a$, coined as the It\^{o}-Stokes drift (ISD) in \cite{Bauer-et-al-JPO-20}, represents through the divergence of the variance tensor, the effects of the small-scale inhomogeneity on the large-scale flow component. This term can be understood as a generalization of the Stokes drift associated to the surface waves velocities.  Note that such a term is added as a corrective advective  term together with its associated vortex force in large scale simulation of ocean dynamics  to take into account the effect of surface waves and Langmuir turbulence \cite{Craik-Leibovich-76,Harcourt-08,McWilliams97}. As shown in \cite{Bauer-et-al-JPO-20}, the LU modelling carries such effects in its own, and generalizes statistical effect of the small-scale inhomogeneity.

Compared to the deterministic material derivative, the stochastic transport operator in \eqref{eq:STO}, involves meaningful terms for large-scale representation of fluid flows. The last term is an inhomogeneous diffusion driven by the variance tensor principal directions. It represents the mixing effect of the small-scale component. Although it can be seen as a matrix generalization of the Boussinesq eddy viscosity assumption, its shape is directly imposed by the noise form. 
 The third right-hand side term corresponds to an energy backscattering from the unresolved scales to the large scales.  The backscattering term corresponds to an energy  source that is exactly compensated by the diffusion term \cite{Resseguier-GAFD-I-17}. 

In particular, for an isochoric flow with $\div (u - u_s) = 0$, one may immediately deduce from \eqref{eq:SRTT} the following transport equation of an extensive scalar:
\begin{equation}
\Df_t q = 0,
\end{equation}
where the energy of such random scalar $q$ can be shown to be globally conserved, through a direct application of It\^{o} formula \cite{Resseguier-GAFD-I-17}:
\begin{equation}
\dif\, \Big( \int_{\sub\mathcal{S}} \frac{1}{2} q^2\, \dif x \Big) = \Big( \underbrace{\frac{1}{2}\int_{\sub\mathcal{S}} q \div (a \nabla q)\, \dif x}_{\text{Energy loss by diffusion}} + \underbrace{\frac{1}{2}\int_{\sub\mathcal{S}} (\nabla q)\tp a \nabla q\, \dif x}_{\text{Energy intake by noise}} \Big)\, \dif t = 0.
\end{equation}
Indeed, this can be interpreted  as a process where the energy brought by the noise is exactly counter-balanced by the energy lost by the diffusion term.

As already mentioned, decomposition \eqref{eq:dX} is written in terms of a It\^{o} stochastic integral. This decomposition could have been written in the form of a Stratonovich integral as well. The calculus associated to this latter integral has the advantage to follow the classical  chain rule. However, the Stratonovich noise is not anymore of zero expectation.  This leads thus to a problematic decomposition with velocity fluctuations of non null ensemble mean. For smooth enough integrands, it is possible to safely move from one form to the other.  For interested readers, more insights on the difference of the two settings and their implications in stochastic oceanic modelling are provided in \cite{Bauer-et-al-JPO-20}.

\subsection{Stochastic Navier-Stokes equation in LU form}

The above SRTT \eqref{eq:SRTT} and Newton's second principle (in a distributional sense) allow us to derive the following  stochastic equations of motions \cite{Memin14, Mikulevicius04}, which for any noise scaling $\varepsilon>0$ parameter and for all points of $\mathcal{S}$ reads, using $\sigma,\; u_s$, and $a$ introduced above: 
\begin{multline}\label{EqNVS}
\dif_t u \, 
+ \, (u- \varepsilon^{2} u_s) \bcdot \nabla u \, \mathrm{d}t \,
+ \, \varepsilon \sigma \mathrm{d}W_t \bcdot \nabla u \, - \, \dfrac{1}{\,2\,} \, \varepsilon^{2} \, \nabla \bcdot (a\nabla u) \, \mathrm{d}t \\
= - \dfrac{1}{\,\rho\, } \nabla(p\,\mathrm{d}t \, +\,  \mathrm{d} p^{\sigma}_t ) \,
+ \, \dfrac{1}{\,R_e\,} \, \Delta (u\,\mathrm{d}t \,
+ \, \varepsilon \sigma \, \mathrm{d}W_t), \qquad 
\end{multline}
with the incompressibility conditions
\begin{equation}
\nabla \bcdot (u- \varepsilon^{2} u_s) = 0 \qquad , \qquad \nabla\bcdot \sigma =0 \; , 
\end{equation}
and associated with Dirichlet boundary condition $u(t,x)=0$ and $\widehat{\sigma}(x,y,t)=0 $ for all $x\in \partial\mathcal{S}$ and $t>0$. The initial condition is denoted by $u(0,x)=u_0(x)$ for all $x\in\mathcal{S}$. As usual, $u(t,x) = \bigl(u_1(t,x) , ... , u_d(t,x)\bigr)$ and $p(t,x)$ stands for the velocity and the pressure of the fluid, respectively. The term $\mathrm{d} p^{\sigma}_t$ corresponds to the Brownian (martingale) part of the pressure. The Ito-Stokes drift $u_s$ is defined as $u_s:= \dfrac{1}{\,2\,} \nabla \bcdot a$ and $\rho$ stands for the fluid density. The dimensioning constant $R_e =U L/\nu$ denotes the Reynolds number, sets from the ratio of the product of characteristic length and velocity scales, $UL$, with the kinematics viscosity $\nu$.  As for the noise scaling parameter, $\epsilon$, it encodes the amplitude scale of the unresolved energy and should converge to zero when  all the flow component are resolved. Meaning thus there is no noise and the system corresponds trivially to the deterministic Navier-Stokes system. In \cite{Resseguier-GAFD-II-17} this factor is defined as the ratio between the turbulent kinetic energy (TKE) and the mean kinetic energy (MKE), multiplied by the ratio between the unresolved scale correlation time $\mathcal{T}_{\sigma}$ and the large-scale advection time. This quantity tends to zero when all the scales are resolved. 

Although the system  corresponds to the Navier-Stokes for zero noise, the convergence  toward weak (strong) solutions of the 3D (2D) deterministic Navier-Stokes, respectively, at the limit of vanishing noise is an important property that should be ideally respected by any stochastic flow representations to ensure physical relevance when all the scales are resolved without almost any uncertainty. This is the main results we aim to prove in this paper.   

First of all, in order to work with a pressure-free system through a divergence-free Leray projection,  we proceed to
the  change of variable $v:= u-\varepsilon^{2} u_s $ in \eqref{EqNVS} to rewrite the system with a classical incompressibility condition on $v$:
\begin{multline}\label{EqNVS2}
d_t v \, 
+ \, v \bcdot \nabla v \, \mathrm{d}t \, 
- \,  \dfrac{1}{\,R_e\,} \Delta v \, \mathrm{d}t \,
+ \, \varepsilon^{2} (v\bcdot \nabla )u_s \, \mathrm{d}t \,
-\,  \dfrac{\,\varepsilon^{2}\, }{\,2\,}  \, \nabla \bcdot (a\nabla v) \, \mathrm{d}t
-\, \dfrac{\,\varepsilon^{4}\,}{\,2\,} \,  \, \nabla \bcdot (a\nabla u_s) \, \mathrm{d}t \, - \, \dfrac{\varepsilon^{2}}{R_e} \Delta u_s \, \mathrm{d}t + \varepsilon^{2} \partial_t u_s \dif t\\
= - \dfrac{1}{\,\rho\, } \nabla(p\,\mathrm{d}t \, +\,  \mathrm{d} p^{\sigma}_t ) \; - \; (\varepsilon \sigma \mathrm{d}W_t \bcdot \nabla)v \; - \; (\varepsilon^{3} \sigma \mathrm{d}W_t \bcdot \nabla) u_s \, + \;  \dfrac{\,\varepsilon\, }{\,R_e\,} \, \Delta ( \sigma \, \mathrm{d}W_t) ,
\end{multline}
 with the incompressibility conditions 
\begin{equation}\label{EqNVS20}
\nabla \bcdot v = 0 \qquad  \qquad \nabla\bcdot \sigma =0 \; , 
\end{equation}
for all points in $\mathcal{S}$ together with Dirichlet boundary conditions $v(t,x)=0$, $\widehat{\sigma}(x,y,t)=0$ for all $x\in \partial\mathcal{S}$, $y\in \mathcal S$ and $t>0$ and the initial condition $v(0,x)=v_0(x):= u_0(x) - \varepsilon^{2}u_{s}(0,x)$ for all $x\in\mathcal{S}$. 
In the following section we specify the spaces on which this system is defined and rewrite it in an equivalent abstract form. 
\section{Preliminaries and main result}
Let us first introduce the functional spaces on which system \eqref{EqNVS} is defined as well as some associated notations used in the following.
\subsection{Definition of the spaces}\label{DefSpace}
Let $\mathcal{V}$ be the space of infinitely differentiable $d$-dimensional vector fields $u$ on $\mathcal{S}$, with compact support strictly contained in $\mathcal{S}$, and satisfying $\nabla \bcdot u = 0$. We denote by  $H$ the closure of $\mathcal{V}$ in $L^{2}(\mathcal{S} , \R^{d})$ and by $V$ the closure of $\mathcal{V}$ in the Sobolev space $H^{1}(\mathcal{S} , \R^{d})$. The space $H$ is  endowed with the $L^{2}(\mathcal{S}, \R^{d})$ inner product. This inner product and its induced norm are writen:
\[
(u,v)_{_H} :=  (u,v)_{L^{2}(\mathcal{S}, \R^{d})} \quad \text{and}\ \quad | u |_{_H} := \|u\|_{L^{2}(\mathcal{S}, \R^{d})} \,.
\]
As for space $V$, thanks to Poincar\'{e} inequality, it is endowed with the $H^{1}_{0}(\mathcal{S}, \R^{d})$ inner product and its associated norm, denoted respectively by
\[
((u,v))_{_V} :=  (\nabla u ,\nabla v)_{L^{2}(\mathcal{S}, \R^{d})} \quad\text{and} \ \quad \| u \|_{_V} := \|\nabla u\|_{L^{2}(\mathcal{S}, \R^{d})}.
\]
We may define then the Gelfand triple $V \subset H \subset V'$ where $V'$ is the dual space of $V$ relative to $H$. We note  $\langle \, \cdot\, , \, \cdot \rangle_{V'\times V}$  the duality pairing between $V'$ and $V$. 

The space of Hilbert-Schmidt operators from the Hilbert space $K_1$ to the Hilbert space $K_2$ is denoted by $\mathcal L_2(K_1,K_2)$  and  $\|\cdot \|_{\mathcal{L}_{2}(K_1,K_2)}$ is its norm.
\subsection{Pressure-free formulation and abstract formulation}
System  \eqref{EqNVS} may be rewritten in an equivalent  simplified pressure-free formulation by using the Leray projection $P :L^{2}(\mathcal{S}, \R^{d})\,\to H$ of $L^{2}(\mathcal{S} , \R^{d})$ onto the space $H$ of divergence-free vectorial functions. Applying Leray's projector to \eqref{EqNVS2}, we obtain 
\begin{multline}\label{EqNVS3}
d_t v \, 
- \,  \dfrac{1}{\,R_e\,} P (\Delta v \, \mathrm{d}t) \, 
+ \, P (v \bcdot \nabla v \, \mathrm{d}t) \\
+ \, P\left(\varepsilon^{2} (v\bcdot \nabla )u_s \, \mathrm{d}t \,
-\,  \dfrac{\,\varepsilon^{2}\, }{\,2\,}  \, \nabla \bcdot (a\nabla v) \, \mathrm{d}t \, 
-\, \dfrac{\,\varepsilon^{4}\,}{\,2\,} \,  \, \nabla \bcdot (a\nabla u_s) \, \mathrm{d}t  \, - \, \dfrac{\varepsilon^{2}}{R_e} \Delta u_s \, \mathrm{d}t + \varepsilon^{2} \partial_t u_s \dif t\right) \\
=  P \left( \dfrac{\,\varepsilon\, }{\,R_e\,} \, \Delta ( \sigma \, \mathrm{d}W_t) \; - \; (\varepsilon \sigma \mathrm{d}W_t \bcdot \nabla)v \; - \; (\varepsilon^{3} \sigma \mathrm{d}W_t \bcdot \nabla) u_s \,  \right) \; .
\end{multline}
This system can finally be rewritten in the following simplified abstract form
\begin{equation}\label{AbstractProblem}
\left\{
    \begin{array}{ll}
        d_t v(t) \; + \; Av(t)  \, \mathrm{d}t \; + B( v(t))\, \mathrm{d}t \; + \; F_{\varepsilon} (v(t)) \, \mathrm{d}t \; = \; G_{\varepsilon} (v(t)) \, \mathrm{d} W_t,  &  \\
        v(0)=v_0. &
    \end{array}
\right.
\end{equation}
The deterministic terms $A$, $B$, $F_{\varepsilon}$ and the stochastic term $G_{\varepsilon}$ are fully described in section \ref{Operators}. We first  present our main results.

\subsection{Main results}
Several kinds of solutions can be defined for stochastic partial differential equations. As for deterministic PDEs, these can be strong, weak or mild (semi-group) solutions. The solutions  are said to be strong in the probabilistic sense when they are constructed for a fixed Wiener process $W$ on a given stochastic basis $(\Omega , \mathcal{F},  (\mathcal{F}_t)_{t\in[0,T]}, \mathbb{P} )$, composed of filtration $(\mathcal{F}_t)_{t\in[0,T]}$ and probability space $(\Omega , \mathcal{F}, \mathbb{P} )$. For stochastic evolution models like the stochastic Navier-Stokes equations  in dimension $3$,  it is more natural to work with weaker solutions, called martingale solutions. In this case, we look for solutions defined as a triplet composed of a stochastic basis, a Wiener process and an adapted process.

More precisely, we say that there is a martingale solution of system (\ref{AbstractProblem}) if there exists a stochastic basis $(\Omega , \mathcal{F},  (\mathcal{F}_t)_{t\in[0,T]}, \mathbb{P} )$, a cylindrical Wiener process $W$ on $L^2(\mathcal S, \R^d)$ and a progressively measurable process $v : [0, T] \times \Omega \to H$, with 
\[
\label{eq-SNS}
v_\varepsilon\in L^{2}\left(\,\Omega \times [0, T] ; V \right) \cap L^{2}\left( \Omega \, , \, {L^{\infty}}([0, T]; H)\right), 
\]
such that $\mathbb{P}-a.e$,  \; $v_\varepsilon$ satisfies for all time $t\in [0,T]$ 
\begin{equation}
\label{Abs-SNS}
v_\varepsilon(t) \; + \; \int_{0}^{t} Av_\varepsilon(s)\, \mathrm{d}s \; + \; \int_{0}^{t} B\bigl(v_\varepsilon(s)\bigr)\,\mathrm{d}s \; + \; \int_{0}^{t} F_\varepsilon \bigl(v_\varepsilon(s)\bigr)\,\mathrm{d}s
\; = \; v_0 \;  + \; \int_{0}^{t} G\bigl(v_\varepsilon(s)\bigr)\,\mathrm{d}W_s,
\end{equation}
where the equality must be understood in the weak sense. Martingale solution are hence the equivalent, in the stochastic setting, of Leray-Hopf weak solution with an additional degree of freedom provided by a complying Wiener process. 
In this work, we will show, for all $\varepsilon>0$, the existence in dimension $2$ or $3$ of a martingale solution for the LU representation of the Navier-Stokes equations for noises associated with a smooth enough  diffusion tensor kernel $\widehat{\sigma}$  in space and time.  In dimension $2$, this solution is unique and strong in the probabilistic sense. 
This result is summarized in the following theorem. 
\begin{theorem}\label{Theorem}  

Let $d=2$ or $3$ and assume that the noise is smooth enough in the sense that its variance tensor and Ito-Stokes drift are such that 
\begin{equation}
\label{hyp_phi}
\displaystyle \qquad \sup_{t\in [0,T]} \|\sigma(t)\|^2_{_{\mathcal L_2(L^2(\mathcal S,{\mathbb R}^d), H^{3}(\mathcal{S}))} }= \quad \sup_{t\in [0,T]} \sum_{k=0}^{\infty} \| \phi_k(t) \|^{2}_{_{H^{3}(\mathcal{S})}} < \infty, 
\end{equation}
\begin{equation}
\label{hyp_a}
u_s \in L^\infty\bigl(0,T;H^{3}(\mathcal{S} , \mathbb{R}^{d})\bigr) \quad ; \quad  \partial_t u_s \in L^\infty(0,T; H) \qquad \text{and} \qquad a \nabla u_s \in L^\infty(0,T;V).
\end{equation}
Then, for all $\varepsilon>0$, equation \eqref{Abs-SNS} admits a martingale solution $v_\varepsilon$. Moreover, for $d=2$, any solution of \eqref{Abs-SNS} is strong in the probabilistic sense, belongs to 
$L^2(\Omega,C([0,T];H))$ and  is  unique. 

Besides, when $\varepsilon \to  0$, for $d=3$, any converging subsequence of $(v_\varepsilon)_{\varepsilon >0}$ converges in law to a solution of the deterministic Navier-Stokes equation in $L^{2}([0,T] \, ; \, H) \, \cap \, C^{0}([0,T] \, ; \, \D(A^{-3/2})\, )$ and there exists such subsequences. For $d=2$, the whole 
sequence converges to the unique solution of the Navier-Stokes equation  in {$L^2(\Omega, L^\infty([0,T];H)\cap L^2([0,T];V))$}.
\end{theorem} 
{\begin{remark}\label{r3.2}
When $d=3$, although the limit equation is deterministic, the limit solution of the Navier-Stokes equations is a random object. This unusual feature is due to the lack of uniqueness for the 3D Navier-Stokes equation. 
For smooth initial data, it is know that there exists a unique strong solution on a  small interval of time $[0,T(u_0))$ and that weak-strong uniqueness holds. As a result, the limit solution obtained in Theorem \ref{Theorem} is not random when restricted to $[0,T(u_0))$.  
\end{remark}}
The proof, which is  thoroughly developed in the following sections, follows the well established approach developed in \cite{Debussche-et-al-2011,Flandoli-Gatarek-95}. It is composed of the following successive steps. 
\begin{itemize}
\item An approximate finite dimensional system obtained from a Galerkin projection is first studied:
\begin{equation*}
\left\{
    \begin{array}{ll}
        d_t v_n(t) \; + \; Av_n(t)  \, \mathrm{d}t \; + B^n_\varepsilon (v_n(t))\, \mathrm{d}t \; + \; F^n_{\varepsilon} v_n(t) \, \mathrm{d}t \; = \; G^n_{\varepsilon} v_n(t) \, \mathrm{d}W_t,  &  \\
        v_n(0)=v_0. & 
    \end{array}
\right.
\end{equation*}
The existence and uniqueness of a solution $v_{\varepsilon,n}$ is proved for all integers $n$.
\item The sequence  $(v_{\varepsilon,n})_n$ is shown to verify energy estimates. 
\item  These estimates enable us to prove that the laws $(\mathcal{L}(v_{\varepsilon,n}) )_n$ are weakly compact (tight) on an appropriate space.
\item We then apply the Skorohod's embedding theorem and change  the stochastic basis. On this basis there exists $(\overline{v}_{\varepsilon,n})_n$ such that $\overline{v}_{{\varepsilon,n}} \overset{\mathcal{L}}{=} v_{\varepsilon,n}$ and by thinning,  the sequence $(\overline{v}_{\varepsilon,n})_n$ converges on an appropriate space almost surely to  $\overline{v_\varepsilon}$, a martingale solution.
\item Finally, in 2D we prove pathwise uniqueness and use a Lemma due to Gyongy and Krylov to deduce that the sequence $(v_{\varepsilon,n})_n$ converges in probability to a strong solution. 
\item Having built  a family $(v_{\varepsilon})_{\varepsilon>0}$ of solutions of \eqref{Abs-SNS}, we prove then that this family converges  to a {- possibly random for $d=3$ - }solution of the deterministic Navier-Stokes equation when $\varepsilon\to 0$. 

\end{itemize}

The condition of Theorem \ref{Theorem} simplifies when the covariance operator does not depend on time or if the ISD is divergent free. In both cases the condition on the temporal derivative of the ISD are not necessary. We note also, that for a spatially homogeneous noise, the variance tensor is constant and the ISD cancels. However this requires  additionally either a periodic domain or the full space.

The assumptions on the noise are anyway non optimal but it is not the purpose of this paper to consider non spatially smooth noise since in practice it is smooth. 
Note that condition \eqref{hyp_phi} is satisfied for instance if we choose $\sigma$ independent on $t$ and equal to $A^{-r}$ with $r$ large enough where $A$ is the Stokes operator defined in the following section. Indeed, in this case $\phi_k= \lambda_k^{-r}e_k$ where $(e_k)_k$ is an orthonormal 
complete system of eigenvectors of $A$ associated to the eigenvalues  $(\lambda_k)_k$ and $\|\phi_k(t) \|^{2}_{_{H^{3}(\mathcal{S})}}= \lambda_k^{3-2r}$. The behavior of the eigenvalues: 
$\lambda_k\sim k^{2/d}$ allows to conclude that \eqref{hyp_phi} follows. Since $u_s= \frac12 \nabla\bcdot a$ and $a$ is defined by \eqref{def_a}, \eqref{hyp_a} holds also for $r$ large enough since 
$\|u_s\|_{H^3(\mathcal S)}\le \sum_{k=0}^{\infty} \| \phi_k(t) \|^{2}_{_{H^{4}(\mathcal{S})}}$. Finally, since $A^{-r}$ is self-adjoint and Hilbert-Schmidt for $r>d/4$, it is associated to a symmetric kernel $\hat \sigma$ which is bounded for $r$ large enough. 

The convergence of the LU representation for limiting vanishing noise warrants that the LU Navier-Stokes equations can be interpreted as a large-scale model of the deterministic Navier-Stokes equation. 

These convergence results open new interesting possibilities for the study of turbulence or for the proposition of new large-scale representations of fluid dynamics. From the theoretical point of view, it might be interesting to explore multiscale versions of the LU representation based on  spatial filtering together with nested noise models. This would generalize classical large eddy models in which the noise would depend on the spatial filtering applied. The coarser the filtering the larger the noise. Energy transfer between scales would then be very interesting to study in this probabilistic setting. Stochastic Kolmogorov-Monin-Horwart equations for energy exchanges across scale could be obtained by this way. From a practical point of view, these convergence results justify the setting of such stochastic models to represent large-scale solution of the Navier-Stokes equations. Compared to the over-diffusive schemes usually employed to that end, these stochastic representations have the advantage to be well suited to uncertainty quantification, ensemble forecasting and data-assimilation \cite{   Bauer-et-al-JPO-20,Bauer-et-al-OM-20, Dufee-et-al-QJRMS-22, Resseguier2020arcme, Resseguier-SIAM-UQ-21, Yang-QJRMS19}, as well as for large-scale flow modelling as used in \cite{Pinier-PRE19} for the characterization of the velocity profile in turbulent boundary layer flows.

\section{Proofs of the main results}\label{MathFramework}

In the following we provide a complete proof of the two results given previously. This section is structured in the following way. We first define the different operators involved in the abstract formulation \eqref{Abs-SNS} and provide estimates for each of them. In a second time we describe the approximate Galerkin system  and give energy estimates. Tightness of the law of the Galerkin approximate solutions is shown in section \ref{Tight} and passage to the limit is then performed. This concludes the proof on the existence (and uniqueness) of martingale (strong) solution in 3D (in 2D), respectively. The convergence toward weak solutions of the deterministic Navier-Stokes equation in the limit of vanishing noise is finally shown in section \ref{PassageLimiteEpsilon}. 

From now on, $C$ denotes a constant which may depend on the domain $\mathcal S$ or the noise characteristics but not on other parameters such as $\varepsilon$ or $N$. 
Also, since we are interested in vanishing noise, we assume that $\varepsilon\in (0,1]$. Our results clearly extend to larger $\varepsilon$. 
 
\subsection{Abstract formulation and Galerkin approximation}\label{Operators}

Let $A : \D(A)\subset H \to H$ be the unbounded linear operator defined by 
\begin{equation}\label{DefA}
Av := -\frac{1}{\,R_e\,} \, P(\Delta v), 
\end{equation}
on the domain $\D(A):=V \cap H^{2}(\mathcal{S} , \mathbb{R}^{d})$. For all $v\in \D(A)$ and $w\in V$, we have
\begin{equation}\label{RelationA}
(Av \, , \, w)_{H} = \frac{1}{R_e} \, ((v \, , \, w ))_{V}, 
\end{equation}
since $P$ is a orthogonal projection onto $H$. Let $a$ be the continuous bilinear form defined by $a(v,w):= \frac{1}{R_e} \int_{\mathcal{S}} \nabla v \bcdot \nabla w \, \mathrm{d}x$ for all $v$ and $w\in V$. For all $v\in V$, we have $Av = a(v\, , \, \cdot)$, in particular, $Av \in V'$ and
\begin{equation}\label{AVVprim}
|Av|_{_{V'}}\leq \frac{1}{R_e} \| v \|_{_V} .
\end{equation}
We also have the equality for all $v$ and $w\in V$, 
\begin{equation}\label{DualiteA}
\langle Av \,, \,  w \, \rangle_{V'\times V} = \dfrac{1}{\,R_e\,} \, ((v \, , \, w))_{_V} \; .
\end{equation}
\noindent
The operator $A : \D(A)\subset H \to H$ is positive self-adjoint with compact resolvent, since the embedding $V\hookrightarrow H$ is compact. By the spectral theorem, we denote by $(\lambda_i)_{i\geq 0}$ the increasing (and unbounded) eigenvalues of $A$ and by $(e_i)_{i\geq 0}$ a corresponding orthonormal Hilbertian basis in $H$ of eigenvectors of $A$.
This basis enables to define new spaces
\[
\D(A^{\alpha}) := \{ v \in H \; : \; \displaystyle \sum_{k=0}^{\infty} \lambda_{k}^{2\alpha} \, | (u\,,\,e_k)_{_H} |^{2} < \infty \, \}
\]
for all $\alpha>0$. We  endow $\D(A^{\alpha})$ with the Hilbertian norm
\[
\| v \|_{_{\D(A^{\alpha})}}:= \left( \displaystyle \sum_{k=0}^{\infty} \lambda_{k}^{2\alpha} \, | (v\,,\,e_k)_{_H} |^{2} \right)^{\alf} \; .
\]
All eigenvectors $e_i$ belong to $\D(A^{\alpha})$ for any $\alpha$. By \eqref{DualiteA}, it can be observed that $V= \D(A^{\alf})$ and $ \| v \|_{V} = \frac{1}{\,R_e\,} \| v \|_{{\D(A^{\alf})}} $ for all $v\in V$. 
Moreover, since the domain $\mathcal S$ is smooth, it can be proved that the norm of the classical Sobolev space  $H^{2\alpha}(\mathcal{S}, \mathbb{R}^{d})$ is equivalent to $\| \bcdot \|_{\D(A^\alpha)}$ .

We have the Gelfand triple $\D(A^{\alpha}) \subset H \subset \D(A^{-\alpha})$ relative to $H$. We denote by $ \langle \cdot \, , \, \cdot \rangle_{\D(A^{-\alpha}) \times \D(A^{\alpha})}$ the duality product. We have, for all $u\in H$ and $v \in \D(A^{\alpha})$,
\[
\langle u \, , \, v \rangle_{\D(A^{-\alpha}) \times \D(A^{\alpha})} \; = \; (u \, , \,  v )_{_H}.
\]
Although incorrect from a strict point of view, we may use the notation $(u,v)_H$  for $u\in \D(A^{-\alpha}),\; v\in \D(A^{\alpha})$.

\medskip
Let $b$ be the trilinear form defined for all $u, v$ and $w\in V$ by  
\[
b(u,v,w) =\bigl(w\,,(u\bcdot \nabla)v\bigr)_H= \displaystyle \int_{\mathcal{S}} \; w(x) \, \left(u(x) \bcdot \nabla \right) v(x)  \, \mathrm{d}x.
\]
By Cauchy-Schwarz,  and H\"older inequality, $b$ satisfies for all $u, v$ and $w\in V$
\[
  | b(u,v,w) | \leq \|u\|_{_{L^{4}(\mathcal{S}, \mathbb{R}^{d})}}  \|w\|_{_{L^{4}(\mathcal{S}, \mathbb{R}^{d})}} \|v\|_{_V}.
\]
When the dimension $d=2$, Gagliardo-Nirenberg inequality and Poincar\'e inequality give for $v\in V$, 
\[
\|v\|_{_{L^{4}(\mathcal{S}, \mathbb{R}^{d})}}\leq C |v|^{\alf}_{_H} \, \|v\|^{\alf}_{_V} \leq C \|v\|_{_V}.
\]
For $d=3$, we have similarly for $v \in V$, 
\[
\|v\|_{_{L^{4}(\mathcal{S}, \mathbb{R}^{d})}}\leq C |v|^{\scriptscriptstyle{1/4}}_{_H} \, \|v\|^{\scriptscriptstyle{3/4}}_{_V} \leq C \|v\|_{_V}.
\]
To conclude, for $d=2$ or $3$, we obtain for all $u,v$ and $w \in V$
\begin{equation}\label{ContinuiteTrilib}
| b(u,v,w) | \leq C \, \|u\|_{_V} \|v\|_{_V} \|w\|_{_V}, 
\end{equation} in particular $b$ is continuous on $V\times V \times V$.  \\

\noindent
Since the Leray projection $P$ is self-adjoint in $L^{2}(\mathcal{S} , \mathbb{R}^{d})$, we have
\begin{equation}\label{DefBPbAstrait}
b(u,v,w) = \int_{\mathcal{S}} P(u\bcdot \nabla v) \, w \; \mathrm{d}x.
\end{equation}
Let $B$ be the bilinear map defined by $B(u,v) :=b(u,v,\, \bcdot)$ for all $u$ and $v\in V$. For the sake of conciseness we write $B(u):=B(u,u)$. Equation \eqref{DefBPbAstrait} provides an explicit expression of $B(v)$ as:
\begin{equation}
\label{DefB2}
    B(v) = P(v \bcdot \nabla v ).
\end{equation}
By \eqref{ContinuiteTrilib}, we have for all $u$ and $v \in V$ that $B(u,v) \in V'$ and 
\begin{equation}\label{BContinuite1}
| B(u,v)|_{V'} \leq C \, \|u\|_{_V} \|v\|_{_V} .
\end{equation}
The bilinear mapping $B : V \times V \to V'$ is hence continuous. \\ 

\noindent
Let $R>0$, for all $u$ and $v\in B_{V}(0,R)$, with $B_{V}(0,R)$ denoting the zero-centered ball of radius $R$ in $V$,  we have
\begin{align*}
 | B(u) - B(v)|_{V'} 
&=  | B(u, u-v) - B(u-v,v)|_{V'} \\
&\leq C\left( \|u\|_{_V} \|u-v\|_{_V} \, + \, \|u-v\|_{_V} \|v\|_{_V}\right) \\
&\leq 2C\, R \, \|u-v\|_{_V} .
\end{align*}
This shows that  $B : V \times V \to V'$ is locally Lipschitz. \\

\noindent
As the function of $V$ are divergence free, for all $u,v$ and $w \in V$ we have
\begin{equation}\label{Egalitesb}
b(u,v,w)=-b(u,w,v) \qquad \text{and} \qquad b(u,v,v)= 0 .
\end{equation}
Consequently, $B : V \times V \to V'$ satisfies for all $u$, $v$ and $w\in V$:
\begin{align}
\langle B(u,v) \, , \, v \rangle_{_{V'\times V}} &=0, \label{PropB1}\\
\langle B(u,v) \, , \, w \rangle_{_{V'\times V}} &= - \langle B(u,w) \, , \, v \rangle_{_{V'\times V}}\; . \nonumber 
\end{align}

\noindent
The operator $B$ can be extended to $H\times H$.  Indeed by the Sobolev embedding $H^{\beta}(\mathcal{S}, \R^d) \hookrightarrow  L^{\infty}(\mathcal{S}, \R^d)$ when $\beta > d/2$, the embedding $\D(A^{\beta/2})\hookrightarrow  L^{\infty}(\mathcal{S}, \R^d)$ is continuous. 
We have hence for all $u,v\in V$ and $w\in \D\left(A^{\frac{\beta+1}{2}}\right)$, 
\begin{align*}
    | b(u,v,w)| = | b(u,w,v)| &\leq \, \| \nabla w \|_{{L^{\infty}(\mathcal{S})}} \, |u |_{{_H}} \, |v |_{{_H}} \\
    &\leq C \, \| \nabla w \|_{_{\D(A^{\beta/2})}} \, |u |_{{_H}} \, |v |_{_H} \\
    &\leq C \, \| w \|_{\D\left(A^{\frac{\beta+1}{2}}\right)} \, |u |_{_H} \, |v |_{_H} \; . 
\end{align*}
Let $\gamma := \frac{\beta +1}{2}>\frac{d+2}{4}$. The tri-linear form $b$ can be extended as a tri-linear operator $b:H\times H \times D(A^{\gamma}) \to \R$. In particular, $B$ can be uniquely extended to an operator still denoted $B : H \times H \to D(A^{-\gamma})$ which verifies for all $u$ and $v \in H$
\begin{equation}\label{EstimB}
   \| B(u,v) \|_{D(A^{-\gamma})} \leq C \,  | u |_{_H} \,  | v |_{_H} .
\end{equation}
As $V\subset H$, we have also for all $u\in H$ and $v\in V$ that
\begin{equation}\label{EstimB1}
   \| B(u,v) \|_{D(A^{-\gamma})} \leq C \,  | u |_{_H} \,  \| v \|_{_V}\; . 
\end{equation}

\medskip

To study the other terms in  equation \eqref{Abs-SNS}, some assumptions on  $a$, $u_s$ and $\sigma$ need to be introduced. As already mentioned, since $\sigma$ is Hilbert-Schmidt, it is compact and there exists $(\phi_k( t) )_k$ an orthogonal basis consisting of the eigenfunctions of the covariance operator for each  $t\in\mathbb{R}$. Also there exists  a set  $\{ (\beta_{t}^{k}) _{t\geq0} \, , \, k \in \mathbb{N} \}$ of independent and identically distributed standard Brownian motions such that  the noise term can be decomposed on this basis as
\begin{equation}\label{DiagoSigma}
\sigma(t) \, \mathrm{d}W_t = \displaystyle \sum_{k=0}^{\infty} \phi_k(t) \, \mathrm{d}\beta^{k}_t
\quad \text{and} \quad a(x,t) \; = \; \displaystyle \sum_{k=0}^{\infty} \phi_k(x,t) \, \phi\transp_k(x,t),
\end{equation}
where each eigenfunction is weighted by its corresponding eigenvalue. It can be noticed that as the noise is divergence free, we have $P\phi_k = \phi_k$ for all $k$.

\noindent
Recall that the assumption on the smoothness of the noise is:
\begin{equation}\label{hyp}\begin{array}{c}
\displaystyle  \sup_{t\in [0,T]} \sum_{k=0}^{\infty} \| \phi_k(t) \|^{2}_{_{H^{3}(\mathcal{S})}} < \infty,  \\
u_s \in L^\infty\bigl(0,T;H^{3}(\mathcal{S} , \mathbb{R}^{d})\bigr), \quad  \partial_t u_s \in L^\infty(0,T; H), \qquad a \nabla u_s \in L^\infty(0,T;V).
\end{array}
\end{equation}
By Sobolev embedding   and by definition of $a$, we deduce 
\begin{equation}\label{hyp2}
\sup_{t\in [0,T]}\sum_{k=0}^{\infty} \| \phi_k(t)\|^{2}_{_{L^{\infty}(\mathcal{S},\mathbb{R}^d)}} < \infty \;, \quad \sup_{t\in [0,T]}\| a(t) \|_V+\| a(t) \|_{_{L^{\infty}(\mathcal{S},\mathbb{R}^d)}} < \infty\;.
\end{equation}
 Also, it can be seen that this implies:
\begin{equation}\label{hyp4}
\sup_{t\in [0,T]}\sum_{k=0}^{\infty} \bigl|B\bigl( \phi_{k}(t) , u_s(t)\bigr) \bigr |^{2}_{_H}  < \infty \; , \quad  \sup_{t\in [0,T]} \| \nabla \bcdot a(t) \|_{_{L^{\infty}(\mathcal{S},\mathbb{R}^d)}} < \infty 
\end{equation}

\medskip

Let us now consider the drift terms in \eqref{EqNVS3}, corresponding to $F_{\epsilon}$ in the abstract problem \eqref{AbstractProblem}. In the following for sake of conciseness, we drop the $\epsilon$ subscript:
\begin{equation}
\label{DefF2}
F(v)  = \varepsilon^{2} B(v,u_s)\,
-\,  \dfrac{\,\varepsilon^{2}\, }{\,2\,}  \, P\bigl(\nabla \bcdot (a\nabla v)\bigr) \, 
-\, \dfrac{\,\varepsilon^{4}\,}{\,2\,} \,  \, P\bigl(\nabla \bcdot (a\nabla u_s)\bigr)   \, - \, {\varepsilon^{2}} A u_s \, + \varepsilon^{2} P\partial_t u_s .
\end{equation} 
By \eqref{AVVprim}, \eqref{BContinuite1}, \eqref{hyp} and \eqref{hyp2}, we have for $v\in V$:
\begin{align}\label{EstimF1}
\|F(v)\|_{_{V'}}&\le C \varepsilon^{2} \,  {\|u_s\|_V}\|v\|_{_V}  + \varepsilon^2  \|a\|_{L^\infty(\mathcal S)} \|v\|_{_V} + \varepsilon^4  \|a\|_{L^\infty(\mathcal S)} \|u_s\|_{_V} +
\frac{\varepsilon^2}{Re} \|u_s\|_{_V} +\varepsilon^2 \|\partial_t u_s\|_{_{V'}}\\
&\le C\varepsilon^{2} \, (\|v\|_{_V} \; +1),
\end{align}

\noindent
and,  for all $u$ and $v\in V$,
\begin{equation}\label{ContinuiteF}
\| F(u) - F(v) \|_{_{V'}} \leq C \varepsilon^{2} \, \| u - v \|_{_V}\;.
\end{equation} 

\noindent

\noindent
We  may also extend continuously $F$ to $H$. Note that for $v\in H, \, w\in \D(A)$:
$$
\biggl(P(\nabla \bcdot \bigl(a\nabla v)\bigr),w\biggr)_{\!\!{_H}}= \bigl(\nabla \bcdot (a\nabla v),w\bigr)_{\!{_H}}=  \bigl( v \, , \, \nabla \bcdot (a \nabla w) \, \bigr)_{\!{_H}}\le C |v|_H  \| w \|_{\D(A)},
$$
so that 
$$
\|P(\nabla \bcdot (a\nabla v))\|_{\D(A^{-1})}\le C |v|_H.
$$
Moreover
$$
|B(v,u_s)|_H\le C  \|\nabla u_s\|_{_{L^{\infty}(\mathcal{S})}} \,  | v |_{_H}\le C | v |_{_H}.
$$
It follows:
\begin{equation}\label{EstimF2}
\| F(v) \|_{\D(A^{-1})} \leq C \varepsilon^{2} \, | v |_{_H} + C \varepsilon^{2},
\quad
\| F(v) - F(w) \|_{\D(A^{-1})} \leq C \varepsilon^{2} \, | v - w |_{_H}.
\end{equation}


\medskip

Let us now examine the martingale term. Let $G$ be defined for all $v\in V$ by 
\begin{align*}
G(v)\, \mathrm{d}W_t \, &= P \left( \dfrac{\,\varepsilon\, }{\,R_e\,} \, \Delta ( \sigma \, \mathrm{d}W_t) \; - \; (\varepsilon \sigma \mathrm{d}W_t \bcdot \nabla)v \; - \; (\varepsilon^{3} \sigma \mathrm{d}W_t \bcdot \nabla) u_s \,  \right) \\
&=\sum_{k=0}^{\infty}  \left( \, \frac{\,\varepsilon\,}{R_e} {P(} \Delta \phi_k {)}\, - \, \varepsilon {P(} (\phi_k \bcdot \nabla) v {)} \, - \, \varepsilon^{3} {P(}(\phi_k \bcdot \nabla) u_s {)}\, \right) \; \mathrm{d}\beta_{t}^{k}\\
&= \, \sum_{k=0}^{\infty} \left( \; -\varepsilon \, A \phi_k \, - \, \varepsilon B(\phi_k,v) \, - \, \varepsilon^{3} B(\phi_k, u_s)  \; \right) \, \mathrm{d}\beta_{t}^{k}  .
\end{align*}
%
In other words,  $G : V \to \mathcal{L}(L^2(\mathcal S,{\mathbb R}^d), H)$ is defined  for $v\in V$ and $\varphi\in L^2(\mathcal S,{\mathbb R}^d)$
\[
G(v) \varphi =  - \varepsilon \, A \sigma\varphi \, - \, \varepsilon  B(\sigma\varphi,v) \, - \, \varepsilon^{3} B(\sigma\varphi,u_s)\; .
\]
Thanks to \eqref{hyp}, we see that for $v\in V$, $G(v)\in \mathcal{L}_2(L^2(\mathcal S,{\mathbb R}^d), H)$ and 
\begin{align}
| G(v) |^{2}_{\mathcal{L}_{2}(L^2(\mathcal S,{\mathbb R}^d) , H)}= \displaystyle \sum_{k=0}^{\infty} | G(v)\phi_k|_{_H}^2 &\leq { C  \sum_{k=0}^{\infty} \left( \varepsilon^{2} \, | A \phi_k |^{2}_{_{H}}  + \varepsilon^{2} \| \phi_k \|^{2}_{_{L^{\infty}(\mathcal{S})}} \, \|v\|^{2}_{_V} + \varepsilon^{6} \| \nabla u_s \|^{2}_{_{L^{\infty}(\Omega)}} | \phi_k|^{2}_{_H} \right)} \nonumber \\
&\leq C\varepsilon^{2} ( \| v \|^{2}_{_V} \; + 1 { + \varepsilon^4} ). \label{EstimG0}
\end{align}

\noindent
 Moreover, for $u$ and $v \in V$,
\[
\| G(u) - G(v) \|^{2}_{\mathcal{L}_{2}(L^2(\mathcal S,{\mathbb R}^d) , H)}   \leq C\varepsilon^{2} \| u- v \|^{2}_{_V},
\]
and $G : V \to \mathcal{L}_{2}(L^2(\mathcal S,{\mathbb R}^d) , H) $ is  Lipschitz continuous.

The mapping $G$ also extends  by continuity to  $G : H \to \mathcal{L}_{2}(L^2(\mathcal S,{\mathbb R}^d) , V')$ which verifies for all $v,\, w\in H$, 
\begin{equation}\label{ContinuiteG2}
| G(v) |^{2}_{\mathcal{L}_{2}(L^2(\mathcal S,{\mathbb R}^d) , V')} 
\leq C\varepsilon^{2} ( \, | v |^{2}_{_H} \, + \,  1 \, ),\quad| G(v) - G(w) |^{2}_{\mathcal{L}_{2}(L^2(\mathcal S,{\mathbb R}^d) , V')} \leq C\varepsilon^{2} \, | v - w |^{2}_{_H} \; .
\end{equation}

%
\medskip

Our construction of solutions  is based on Galerkin approximation. We use the Hilbertian basis $(e_i)_{i\geq 0}$ of $H$ consisting of eigenvectors of $A$, and define, for $n\in \mathbb N$, the orthogonal  projections $P_n : \D(A^{-\gamma}) \to H_n := Span (e_0 , ... \, ,  e_n)$:
\[
P_n u := \displaystyle \sum_{i=0}^{n} \langle u \, , \, e_i\rangle_{{\D(A^{-\gamma})\times \D(A^{\gamma})}} e_i,  \; u\in  \D(A^{-\gamma}),\; \gamma\in \R. 
\]
Clearly:
\begin{equation}\label{InegProjection}
\|P_n u\|_{\D(A^{-\gamma})} \leq  \| u\|_{\D(A^{-\gamma})}, \quad \text{$u\in \D(A^{-\gamma})$ } \, .
\end{equation}



\noindent
By dominated convergence, we have:
\begin{equation}
\|P_nv \, - \, v  \|_{D(A^{\gamma})} \underset{n\to+\infty}{\longrightarrow} 0,\quad  v \in D(A^{\gamma}), \quad \gamma\in \R. \label{ConvergenceProjection} \\
\end{equation}

\noindent
Also for $\alpha < \beta$ and all $n\in \N$, 
\begin{equation}\label{ProjProp1}
\| P_n v \|_{\D(A^{\beta})}\leq \, \lambda_{n+1}^{\beta - \alpha} \, \| P_n v \|_{\D(A^{\alpha})},
\quad 
\|\left( I - P_n \right) v \|_{\D(A^{^{\alpha})}} \leq \, \lambda_{n+1}^{\alpha - \beta} \, \| \left( I - P_n \right) v \|_{\D(A^{\beta})} \; .
\end{equation}


 We now  introduce the projected operators:
\[
B^{n}:=P_n B \qquad F^n =P_n F \qquad G^n = P_n G \; .
\]
It can be seen that $B^{n}$ is locally Lipschitz and  $F^{n}$ and $G^{n}$ are globally Lipschitz on $H_n$, thus 
\begin{equation}\label{EquationApprochee}
\left\{
    \begin{array}{ll}
        d_t v_{n} (t) \; + \; Av_{n}(t)  \, \mathrm{d}t \; + B^{n}\bigl(v_{n}(t)\bigr)\, \mathrm{d}t \; + \; F^{n}\bigl(v_{n}(t)\bigr)\, \mathrm{d}t \; = \; G^{n}\bigl(v_n(t)\bigr)\, \mathrm{d}W_t , &  \\
        v_{n}(0)=P_n(v_0) , & 
    \end{array}
\right.
\end{equation}
admits also a local solution $v_{n} \in C^{0}([0,t_n] \, , \, H_n)$. Note that $t_n$ is {\it a priori} random, it can be chosen as a stopping time. We prove in section \ref{PreuveEstimees} that for all $n\in\mathbb{N}$, $v_n$ verifies for all $p\geq 2$:
\begin{equation}\label{Estimees}
\mathbb{E}\left[ \, \sup_{0\leq s \leq t_{n}} \,  |v_n(s)  |^{p}_{_H} \, \right] < C, \qquad \mathbb{E}\left[ \, \displaystyle \int_{0}^{t_{n}} \| v_n(s) \|^{2}_{_V}  \, \mathrm{d}s \, \right] < C \; .
\end{equation}
These energy estimates enable us  to deduce by standard arguments  that the solutions $v_n$ are global in time, in particular $v_{n} \in C^{0}([0,T]\, , \,H_n) $, $t_n=T$  and for all $t \in [0,T]$,  we have,
\begin{equation}\label{SolApprox}
    v_{n}(t) - P_n(v_0) \, + \, \displaystyle \int_{0}^{t} \left[\, Av_{n}(s) +B^{n}\bigl(v_{n}(s)\bigr) +F^{n}\bigl(v_{n}(s)\bigr) \, \right] \, \mathrm{d}s \; = \; \displaystyle \int_{0}^{t} G^{n}\bigl(v_{n}(s)\bigr) \, \mathrm{d}W_s \; . 
\end{equation}

\subsection{Energy estimates}\label{PreuveEstimees}
Apply It\^{o} formula with $F(x)= |x|^{p}_{_H}$ for some $p\geq 2$ to the semi martingale $v_n$:
\begin{multline}\label{EqEstimee1}
d_t |v_{n}(t)|_{_H}^{p} = p |v_{n}(t)|_{_H}^{p-2} \, \biggl( v_n(t) \, , \, G^{n}\bigl(v_n(t)\bigr) \mathrm{d}W_t \biggr)_{\!\!{_H}} \\
- p |v_{n}(t)|_{_H}^{p-2} \, \biggl( v_{n}(t) \, , \, Av_{n}(t) + B^{n}\bigl(v_{n}(t)\bigr) + F^{n}\bigl(v_{n}(t)\bigr) \, \biggr)_{\!\!{_H}} \, \mathrm{d}t \\
+ \dfrac{p(p-2)}{2} \left| G^{n}\bigl(v_n(t)\bigr)^*   v_n(t) \right|^{2}_{_{L^2(\mathcal S,{\mathbb R}^d)}} |v_{n}(t)|_{_H}^{p-4} \mathrm{d}t + \dfrac{p}{2} \|G^{n}(v_{n})(t)\|_{\mathcal{L}_{2}(L^2(\mathcal S,{\mathbb R}^d) , H)} ^{2} |v_{n}(t)|_{_H}^{p-2} \mathrm{d}t\,.
\end{multline}
As $P_n : H \to H_n$ is an orthogonal projection for the inner product of $H$, by \eqref{RelationA} and \eqref{PropB1} we have 
\[
\bigl( v_n(t) \, , \, A v_n(t) \bigr)_{H} = \dfrac{1}{\, R_e\,} \, \| v_n(t) \|^{2}_{_V} \quad \text{and} \quad \biggl( v_n(t) \, , \, B^{n} \bigl(v_n(t)\bigr) \biggl)_{H}=\biggl( v_n(t) \, , \, B \bigl(v_n(t)\bigr) \biggr)_{H} =0 .
\]
We have also
\begin{multline*}
\biggl( v_{n}(t) \, , \,  F^{n}\bigl(v_{n}(t)\bigr) \biggr)_{_H} = \varepsilon^{2} b(v_{n}(t), u_s  , \, v_{n}(t) ) \,-\,  \dfrac{\,\varepsilon^{2}\, }{\,2\,}  \, \biggl( v_n(t) \, , \,  \nabla \bcdot \bigl(a\nabla v_n(t)\bigr) \biggr)_{_H} \\
-\, \dfrac{\,\varepsilon^{4}\,}{\,2\,} \, \bigl( v_{n}(t) \, , \, \nabla \bcdot (a\nabla u_s)\bigr)_{_H}   \, + \, \varepsilon^{2}\bigl( Au_s \, , \, v_{n}(t) \bigr)_{_H} \, + \, \varepsilon^{2} \bigl( \partial_t u_s \, , \, v_n(t) \bigr)_{_H} \\
 :=  \, F_{1}^{n} \, + \, F_{2}^{n} \, + \, F_{3}^{n} \, + \, F_{4}^{n}\, + \, F_{5}^{n}\,. \qquad \ \quad \ \quad \ \qquad \ \qquad \ \qquad
\end{multline*}
For $F_{3}^{n}$, we apply Cauchy-Schwarz  inequality to get, thanks to  \eqref{hyp},
\[
|F_{3}^{n} | \leq \dfrac{\varepsilon^{4}}{2} \, |v_{n}(t)|_{_H}  \, \| a \nabla u_s \|_{_V} \leq C \varepsilon^{4} |v_{n}(t)|^{2}_{_H} + C\varepsilon^{4} \| a \nabla u_s \|^{2}_{_V}
\le  C \varepsilon^{4}( |v_{n}(t)|^{2}_{_H} +1).
\]
We treat the three last terms similarly and get:
\[
| F_{1}^{n} + F_{3}^{n} + F_{4}^{n}  + F_{5}^{n} | \leq \, C \, \varepsilon^{2}\bigl(  \; |v_{n}(t)  |^{2}_{_H} +1\bigr).
\]
By definition of $a$, we have 
\begin{align*}
F_{2}^{n} 
&= \dfrac{\varepsilon^{2}}{2} \bigl(  a \nabla v_{n}(t) \, , \, \nabla v_{n}(t) \bigr)_{_H}\\
&=  \dfrac{\varepsilon^{2}}{2} \displaystyle \sum_{i,j=1}^{d} \int_{ \mathcal{S}} \, \sum_{k=0}^{\infty} \phi_{k}^{i} (x) \, \phi_{k}^{j} (x) \, \partial_j v_{n}(t,x) \partial_i v_{n}(t,x)  \mathrm{d}x \\
&=  \dfrac{\varepsilon^{2}}{2} \displaystyle \sum_{k=0}^{\infty}  | (\phi_k \bcdot \nabla) v_{n}(t) |^{2}_{L^2(\mathcal S,{\mathbb R}^d)}.
\end{align*}
Setting $\psi_k := -\varepsilon^{2} B(\phi_k,  u_s) - A \phi_k$, we can  write by \eqref{PropB1} and \eqref{hyp}:
\begin{align*}
\frac{1}{\,2\,} \; \|G^{n}&v_{n}(t) \|_{\mathcal{L}_{2}\bigl(L^2(\mathcal S,{\mathbb R}^d)\bigr) , H)} ^{2} 
=  \frac{\,\varepsilon^{2}\,}{2}  \; \displaystyle \sum_{k=0}^{\infty} | \psi_k(t) - B(\phi_k, v_{n})(t) |^{2}_{_H} \\
&\leq  \frac{\,\varepsilon^{2}\,}{2}  \; \displaystyle \sum_{k=0}^{\infty} |B(\phi_k, v_{n})(t) |^{2}_{_H} + |\psi_k(t)|^{2}_{_H}  + 2 \, \left| \left(B(\psi_k(t)  , \phi_k ), v_n(t)\right)_{_H} \right| \\
&\leq \frac{\,\varepsilon^{2}\,}{2}  \;  \left( \displaystyle \sum_{k=0}^{\infty} |(\phi_k \bcdot \nabla) v_{n}(t)|^{2}_{L^2(\mathcal S,{\mathbb R}^d)} + |\psi_k(t)|^{2}_{_H} +  |B(\psi_k(t)  , \phi_k )|^{2}_{H} \right) + 2 \varepsilon^{2} \,  |v_{n}(t)|^{2}_{_H} \\
&\leq \frac{\,\varepsilon^{2}\,}{2} \displaystyle \sum_{k=0}^{\infty} |(\phi_k \bcdot \nabla) v_{n}(t)|^{2}_{L^2(\mathcal S,{\mathbb R}^d)} \; + \; C \varepsilon^{2}\bigl(  |v_{n}(t)|^{2}_{_H} +1\bigr) 
\end{align*}
and the first term is exactly  $F_{2}^{n}$.

We now estimate the second last term in \eqref{EqEstimee1}. Write, thanks to Parseval and \eqref{hyp}:

\begin{align}
\bigl|G^{n}\bigl(v_n(t)\bigr)^*   v_n(t) \bigr|^{2}_{_H}&=  \displaystyle \sum_{k=0}^{\infty} \bigl(v_n(t),  \varepsilon \, A \phi_k + \, \varepsilon B(\phi_k,v_n) \, + \, \varepsilon^{3} B(\phi_k, u_s)\bigr)_H^2\\
&=  \displaystyle \sum_{k=0}^{\infty} \bigl(v_n(t),  \varepsilon \, A \phi_k +  \, \varepsilon^{3} B(\phi_k, u_s)\bigr)_H^2\\
&\le
 C\varepsilon^2 |v_n|_H^2.
 \label{IntermBDG}
\end{align}
It follows:
 \begin{equation*}\label{InterEstim}
|v_n(t)|^{p-4}_{_H} \;  \bigl(G^{n}v_n(t) ,  v_n(t) \bigr)^{2}_{_H} \le
 C\varepsilon^2 |v_n|_H^{p-2} \; \leq \; \varepsilon^{2} |v_{n}(t)|^{p}_{_H} \; + \; C\, \varepsilon^{2}.
\end{equation*}
Gathering the previous estimates yields
\begin{align}\label{EqEstimee5}
d_t |v_{n}(t)|_{_H}^{p} + \dfrac{p}{\,R_e\,} |v_{n}(t)|_{_H}^{p-2} \|v_{n}(t)\|_{_V}^{2} 
&\leq  p |v_{n}(t)|_{_H}^{p-2} \, \biggl( v_n(t) \, , \, G^{n}\bigl(v_n(t)\bigr) \mathrm{d}W_t \biggr)_{_H} \\
&+ \;  C \, \varepsilon^{2}  \, |v_{n}(t)|^{p}_{_H} \; + \; C\, \varepsilon^{2}
\end{align}
with $C>0$  not depending on $\varepsilon$ or $n$. Integrating in time gives for $t\in [0,T]$:
\begin{multline}\label{EqEstimee2}
|v_{n}(t)|_{_H}^{p} + \dfrac{p}{\,R_e\,} \int_{0}^{t} |v_{n}(r)|_{_H}^{p-2} \|v_{n}(r)\|_{_V}^{2} \, \mathrm{d}r \; \leq  \; |v_{0}|_{_H}^{p} \; + \; C\, \varepsilon^{2}\\
+ p \int_{0}^{t} |v_{n}(r)|_{_H}^{p-2} \, \biggl( v_n(r) \, , \, G^{n}\bigl(v_n(r)\bigr) \mathrm{d}W_r \biggr)_{\!\!\!{_H}} \; + \;  C \, \varepsilon^{2}  \, \int_{0}^{t} |v_{n}(r)|^{p}_{_H} \, \mathrm{d}r .
\end{multline}
Since $\displaystyle \biggl(  \int_{0}^{t} \biggl(v_n(r) \, , \, G^{n}\bigl(v_n(r)\bigr) \mathrm{d}W_r \biggr)_{\!\!\!{_H}}\biggr)_{\!\!t}$ is a martingale, we obtain
\begin{multline}\label{EqEstimee3}
\mathbb{E}\left[ \,  |v_{n}(t)|_{_H}^{p} + \dfrac{p}{\,R_e\,} \int_{0}^{t} |v_{n}(r)|_{_H}^{p-2} \|v_{n}(r)\|_{_V}^{2} \, \mathrm{d}r \, \right]  \leq \mathbb{E}\left[ \,  |v_{0}|_{_H}^{p} \right]  
+ C \, \varepsilon^{2}  \, \mathbb{E}\left[ \, \int_{0}^{t} |v_{n}(r)|_{_H}^{p} \, \mathrm{d}r \, \right] \; + \; C\, \varepsilon^{2} .
\end{multline}
Taking $p=2$, this estimate implies the second bound of \eqref{Estimees}.

From Gronwall's lemma, there exists a constant $C>0$ such that
\begin{equation}\label{EstimeeEnergie1}
\mathbb{E}\left[ \,  |v_{n}(t)|_{_H}^{p} \, \right] \leq C\varepsilon^{2} ,
\end{equation}
for all $t\in [0,T]$. Then apply the Burkholder-Davis-Gundy (BDG) inequality (Proposition \ref{BDG}) and use  \eqref{IntermBDG}  to get
\begin{align*}
&\displaystyle \mathbb{E}\left[ \sup_{0\leq t \leq T} \left| \int_{0}^{t} |v_{n}(r)|_{_H}^{p-2} \; \biggl( v_n(r) \, , \, G^{n}\bigl(v_n(r)\bigr) \mathrm{d}W_r \biggr)_{_H} \right| \, \right] \\
&\leq C \,  \mathbb{E}\left[ \, \left( \int_{0}^{T} |v_{n}(r)|_{_H}^{2(p-2)} \; \left|G^{n}v_n(r)^* v_n(r)  \right|^{2}_{_{L^2(\mathcal S,{\mathbb R}^d)}} \, \mathrm{d}r \, \right)^{\alf} \right] \\
&\leq C \, \varepsilon^{2}  \; \mathbb{E}\left[ \, \left( \int_{0}^{T} |v_{n}(r)|_{_H}^{2(p-1)} \, \mathrm{d}r \, \right)^{\alf} \right] \\
&\leq C \,\varepsilon^{2}  \;  \mathbb{E}\left[ \,  \sup_{0\leq t \leq T} |v_{n}(t)|_{_H}^{p-1} \, \right] \\
&\leq \dfrac{1}{2} \, \mathbb{E}\left[ \sup_{0\leq t \leq T} |v_{n}(t)|_{_H}^{p} \right] + C\varepsilon^{2}.
\end{align*}
By \eqref{EqEstimee2} we have hence 
$$
\displaystyle \dfrac{1}{2}\;\mathbb{E}\left[ \sup_{0\leq t \leq T} |v_{n}(t)|_{_H}^{p} \right] 
\leq   \mathbb{E}\left[ \,  |v_{0}|_{_H}^{p} \right] + C \, \varepsilon^{2}  \,  \int_{0}^{T} \mathbb{E}\left[|v_{n}(t)|^{p}_{_H} \right] \, \mathrm{d}t 
+ \; C\, \varepsilon^{2} .
$$
By \eqref{EstimeeEnergie1}, we obtain also the first estimate \eqref{Estimees} for all integers $p\geq 2$.

\subsection{Tightness of the laws of the approximation}\label{Tight}
We have for all integers $n$, 
\begin{align}\label{SolApprox2}
    v_{n}(t) 
    &= P_n(v_0) \, - \, \displaystyle \int_{0}^{t} \biggl( Av_{n}(r)  + B^{n}(v_{n}(r)\bigr) + F^{n} \bigl(v_{n}(r)\bigr) \, \biggr) \, \mathrm{d}r \, + \, \displaystyle \int_{0}^{t} G^{n}\bigl(v_{n}(r)\bigr) \, \mathrm{d}W_r \\
    v_{n}(t) &:= J_{n}^{1} \, + \, J_{n}^{2} \, + \, J_{n}^{3} \, + \, J_{n}^{4} \, + \, J_{n}^{5}.
\end{align}
We prove in the following that the family $\left( \mathcal{L}(v_n)\right)_n$ are tight in $L^{2}([0,T] ; H) \, \cap \, C^{0}([0,T] ; \D(A^{-3/2})\, )$. We choose $\gamma$ defined in \eqref{EstimB} such that $\gamma \in ( \frac{5}{4} , \frac{3}{2})$. \\

\noindent
\underline{\textbf{The laws $\bigl(\mathcal{L}(v_n)\bigr)_n$ are tight in $L^{2}([0,T] \, ; \, H)$.}} \\
Let $\alpha \in (0, \frac{1}{2})$ be fixed. By  lemma \ref{Compact1}, the embedding
\begin{center}
$L^{2}([0,T] ; V) \, \cap \, W^{\alpha,2}\bigl([0,T] ; \D(A^{-\gamma}) \,\bigr) \hookrightarrow L^{2}([0,T] ; H)$ is compact.
\end{center}
By \eqref{Estimees}, the laws $\bigl(\mathcal{L}(v_n)\bigr)_n$ are bounded in probability in $L^{2}([0,T] ; V)$. We now prove that these laws are bounded in probability in $W^{\alpha,2}\bigl([0,T] ; \D(A^{-\gamma})\,\bigr)$. As $u_0 \in H$ and $P_n$ is an orthogonal projection in $H$, we have
\[
\displaystyle \mathbb{E}\left[ \, | J_n^{1} |^{2}_{_H} \, \right] \; \leq C \; .  
\]
We have by the definition of the norm of $W^{1,2}(0,T ; V')$, \eqref{AVVprim} and \eqref{Estimees}:
\[
\mathbb{E}\left[ \, \| J_n^{2} \|^{2}_{W^{1,2}([0,T] \,  ; \, V')} \, \right]\; \leq C \; \mathbb{E}\left[ \, \int_{0}^{T} \|A v_n(r)  \|^{2}_{_{V'}} \, \mathrm{d}r \, \right] 
\leq C \; \mathbb{E}\left[ \, \int_{0}^{T} \| v_n(r)  \|^{2}_{_{V}} \, \mathrm{d}r \, \right] \leq C  .
\]
Furthermore, we have 
\[
\| J_n^{3} \|^{2}_{W^{1,2}\bigl(0,T ; \D(A^{-\gamma})\bigr)} \; \leq \; C \,  \int_{0}^{T} \| B \bigl(v_n(r)\bigr) \|^{2}_{_{\D(A^{-\gamma})}} \,\mathrm{d}r \; , 
\]
and
\[
\| J_n^{4} \|^{2}_{W^{1,2}(0,T ; V')} \, \leq C \,  \int_{0}^{T} \| F\bigl( v_n(r)\bigr) \|^{2}_{_{V'}} \,\mathrm{d}r \, .
\]
This implies, by \eqref{EstimB1}, \eqref{EstimF1}, \eqref{Estimees} and Cauchy-Schwarz,
\begin{align*}
&\mathbb{E}\left[ \| J_n^{3} \|^{1/2}_{W^{1,2}([0,T] \, ; \, \D(A^{-\gamma})\, )} \right] \leq \; C  \, \mathbb{E}\left(\sup_{0\leq r \leq T} | v_n(r) |_{_{H}} \; \left(\int_{0}^{T}  | v_n(r) |^{2}_{_{V}}\, {\mathrm{d}r} \right)^{1/2}\right)\, \le C,\\
& \mathbb{E}\left[ \, \| J_n^{4} \|^{2}_{_{W^{1,2}([0,T] \, ; \,V')}} \, \right]\; \leq 2 \, C \, \varepsilon^{2} \, \mathbb{E}\left[ \, \int_{0}^{T} \| v_{n}(r)  \|^{2}_{_V} \,\mathrm{d}r \, \right]\; + \; 2 \,C \, \varepsilon^{2} \; \leq \; 2\,C\varepsilon^{2}.
\end{align*}

Concerning $J_n^{5}$, since $\alpha< 1/2$, we may apply lemma \ref{Lemme21}, there exists a constant $C_{\alpha}>0$ independent of $G$ , $v_n$ and $n$ such that
\[
\mathbb{E}\left[ \, \| J_n^{5} \|^{2}_{W^{\alpha,2}([0,T] \, ; \,  H)} \, \right] \;\leq C_{\alpha} \; \mathbb{E}\left[ \, \int_{0}^{T} \| G^{n}(v_n)(r) \|^{2}_{\mathcal{L}_{2}(L^2(\mathcal S,{\mathbb R}^d) , H)} \, \mathrm{d}r \; \right] .
\]
As $P_n : H \to H_n$ is an orthogonal projection and by \eqref{EstimG0}, we get
\begin{align*}
\mathbb{E}\left[ \, \| J_n^{5} \|^{2}_{W^{\alpha,2}([0,T] \, ; \, H)} \, \right] \;
&\leq C \; \mathbb{E}\left[ \, \int_{0}^{T} \| G(v_n)(r) \|^{2}_{\mathcal{L}_{2}(L^2(\mathcal S,{\mathbb R}^d) , H)} \, \mathrm{d}r \; \right] \; \\
&\leq \; C \varepsilon^{2} \, \mathbb{E}\left[ \, \int_{0}^{T} \|v_n(r) \|^{2}_{_V} \mathrm{d}r \; + \; 1\; \right] .
\end{align*}
Estimate \eqref{Estimees} implies
\[
\mathbb{E}\left[ \, \| J_n^{5} \|^{2}_{W^{\alpha,2}([0,T] \,  ; \, H)} \, \right] \;\leq C \,\varepsilon^{2} \; . 
\]
As $H \subset V' \subset \D(A^{-\gamma})$ and $\alpha < \frac{1}{2}$, we obtain by \eqref{SolApprox2} 
\[ 
\mathbb{E} \left[ \, \| v_n \|^{1/2}_{W^{\alpha,2}([0,T] \,  ; \, \D(A^{-\gamma}) \, )} \, \right] \, \leq \, C\, \varepsilon^{2} 
\]
and by Lemma \ref{Compact1}, the laws of $(v_n)_n$ are thus tight in $L^{2}([0,T] \, ; \, H)$. \\

\noindent
\underline{\textbf{The law $\bigl(\mathcal{L}(v_n)\bigr)_n$ are tight in $C^{0}\bigl([0,T] \, ; \, \D(A^{-3/2})\, \bigr)$. }} \\
Let $\alpha \in (\frac{1}{3} \, ; \, \frac{1}{2})$  be fixed. Since $\D(A^{-\gamma}) \hookrightarrow \D(A^{-3/2})$ is compact and $3 \alpha >1$, by lemma \ref{Compact2}, we infer the compact embeddings
\begin{center}
$W^{1,2}\bigl([0,T] \, ; \, \D(A^{-\gamma})\bigr) \hookrightarrow C^{0}\bigl([0,T] \, ; \, \D(A^{-3/2})\bigr)$ \\
$W^{\alpha,3}\bigl([0,T] \, ; \, \D(A^{-\gamma})\bigr) \hookrightarrow C^{0}\bigl([0,T] \, ; \, \D(A^{-3/2})\bigr)$.
\end{center}

\noindent
By the previous proof, we have directly 
\begin{equation}\label{Tight1}
\mathbb{E}\left[ \, \left\| v_n(t) - \displaystyle \int_{0}^{t} G^{n}(v_{n}(r)) \, \mathrm{d}W_r \, \right\|^{1/2}_{W^{1,2}([0,T] \, ; \, \D(A^{-\gamma}) )} \, \right] \leq C\,  \; .
\end{equation}
Arguing as above for $J_{n}^{5}$, we have
\[
\mathbb{E}\left[ \, \left\| \displaystyle \int_{0}^{t} G^{n}(v_{n}(r)) \, \mathrm{d}W_r \,  \right\|^{3}_{W^{\alpha,3}([0,T] \, ; \, \D(A^{-\gamma}))} \, \right] \; \leq C_{\alpha} \; \mathbb{E}\left[ \, \int_{0}^{T} \| G^{n}(v_n)(r) \|^{3}_{\mathcal{L}_{2}(L^2(\mathcal S,{\mathbb R}^d) , \D(A^{-\gamma}))} \, \mathrm{d}r \; \right]
\]
with $C_{\alpha}>0$ independent of $G$ , $v_n$ and $n$. As $V'\subset \D(A^{-\gamma})$, we have by \eqref{InegProjection}  and \eqref{ContinuiteG2} 
\begin{align}
\mathbb{E}\left[ \, \left\| \displaystyle \int_{0}^{t} G^{n}(v_{n}(r)) \, \mathrm{d}W_r \,  \right\|^{3}_{W^{\alpha,3}([0,T] \, ; \, \D(A^{-\gamma}))} \, \right] 
& \leq C_{\alpha} \; \mathbb{E}\left[ \, \int_{0}^{T} \| G(v_n)(r) \|^{3}_{\mathcal{L}_{2}(L^2(\mathcal S,{\mathbb R}^d) , V')} \, \mathrm{d}r \; \right] \nonumber\\
& \leq C_{\alpha} \,\varepsilon^{3} \;  \mathbb{E}\left[ \, 1 + \sup_{0\leq r\leq T} | v_n(r) |^{3}_{_H} \, \right] \nonumber \\
&\leq \, C_{\alpha}\,\varepsilon^{3}  \; . \label{Tight2}
\end{align}
This implies by \eqref{Tight1} , \eqref{Tight2} and the compact embeddings that the laws of $(v_n)_n$ are hence tight in $C^{0}\bigl([0,T] \, ; \, \D(A^{-3/2})\bigr)$.

\subsection{Passage to the limit}\label{PassageLimite}

Let $\beta>d/2$, $\alpha<1/2$, the Wiener process lives in $C^\alpha([0,T],D(A^{-\beta}))$. Define $W_n=W$, then the family of the laws of the couple $(v_n,W_n)$ is tight in 
$\left(L^{2}([0,T] ; H) \cap C^{0}([0,T] ; \D(A^{-3/2}))\right)\times C^\alpha([0,T],D(A^{-\beta}))$.
By the Skorohod's embedding theorem, there exists a stochastic basis $(\overline\Omega , \overline{\mathcal{F}} , (\overline{\mathcal{F}}_t)_t , \overline{\mathbb{P}} )$ with  $L^{2}([0,T] ; H) \cap C^{0}([0,T] ; \D(A^{-3/2}))$-valued random variables $\overline{v}_n$ for $n\geq 1$ and $\overline{v}$ such that $\overline{v}_n$ has the same law as $v_{n}$ on $L^{2}([0,T] ; H) \cap C^{0}([0,T] ; \D(A^{-3/2}))$ and  cylindrical Wiener processes on $L^2(\mathcal S,\mathbb{R}^d)$ $\overline{W}^{n}$ for $n \geq 1$ together with $\overline{W}$ such that (by thinning the sequences)
\begin{align}
&\overline{v}_n \to \overline{v} \; \text{in} \;  L^{2}([0,T] \, ; \, H) \cap C^{0}([0,T] \, ; \, \D(A^{-3/2})), \qquad  \overline{\mathbb P} \;  \text{a.s} , \label{Convergence} \\
&\overline{W}^{n} \to \overline{W} \; \text{in} \; C^\alpha([0,T],D(A^{-\beta})), \qquad \overline{\mathbb{P}} \;  \text{a.s} \; . \label{ConvergenceSto}
 \end{align}
For all integers $n$, $\overline{v}_n$ verifies 
\begin{equation}\label{EqApprox2}
\overline{v}_{n}(t) - P_n(v_0) \, + \, \displaystyle \int_{0}^{t} \bigl[ A\overline{v}_{n}(r) +B^{n} \overline{v}_{n}(r)  \, + \, F^{n} \overline{v}_{n}(r) \bigr] \, \mathrm{d}r \; 
= \; \displaystyle \int_{0}^{t} G^{n}\bigl(\overline{v}_{n}(r)\bigr) \, \mathrm{d}\overline{W}^{n}_r \; . 
\end{equation}
Since $\overline{v}_n \overset{\cal{L}}{=} v_n$, by the estimates (\ref{Estimees}), we have, for $p\ge 2$,
\begin{align}
    \sup_{n\in\mathbb{N}} \, \overline{\mathbb{E}} \left[ \, \sup_{r \in [0,T]} |\overline{v}_n(r)  |^{p}_{_H} \, \right] \leq C,
    \qquad \sup_{n\in\mathbb{N}} \, \overline{\mathbb{E}} \left[ \,  \displaystyle \int_{0}^{T} \| \overline{v}_n(r) \|^{2}_{_V}  \, \mathrm{d}r \; \right] \leq C \; . \label{Estimevn1}
\end{align}
Therefore, $(\overline{v}_n)_n$ converges weakly to $\overline{v}$ in $L^{2}(\overline{\Omega} ; L^{2}([0,T] , V))$ and weakly star in $L^{2}(\overline{\Omega} ; L^{\infty}([0,T] ; H))$. In particular, $\overline{v}$ belongs to
$ L^{2}(\overline{\Omega} ; L^{2}([0,T] \, , \, V)) \cap L^{2}(\overline{\Omega} ; L^{\infty}([0,T] \, , \, H))$ and verifies
\begin{align}
    \sup_{n\in\mathbb{N}} \; \overline{\mathbb{E}} \left[ \, \sup_{r \in [0,T]} |\overline{v}(r)  |^{2}_{_H} \, \right] &\leq C, 
    \qquad \sup_{n\in\mathbb{N}} \; \overline{\mathbb{E}} \left[ \,  \displaystyle \int_{0}^{T} \| \overline{v}(r) \|^{2}_{_V}  \, \mathrm{d}r \; \right] &\leq C \; .\label{EStimV1}
\end{align} 

\noindent
We pass now to the limit $[n\to\infty]$ in (\ref{EqApprox2}). To that end, we consider the variational problem associated to (\ref{EqApprox2}):  for all $n\in\mathbb{N}$, all $t\in [0,T]$ and all $z\in  \D(A^{3/2})$, 
\begin{equation*}
\bigl( \overline{v}_{n}(t) - P_n(v_0) \, , \, z\bigr)_{\!\!{_H}} \, + \, \displaystyle \int_{0}^{t} \bigl(A\overline{v}_{n}(r) \, , \, z \bigr)_{\!\!{_H}} \mathrm{d}r \, + \, \displaystyle \int_{0}^{t} \bigl( B^{n} \overline{v}_{n}(r) \, , \, z \bigr)_{\!\!{_H}} \mathrm{d}r 
\end{equation*}
\begin{equation}\label{PbVar1}
+ \, \displaystyle \int_{0}^{t} \bigl( F^{n} \overline{v}_{n}(r) \, , \, z\bigr)_{\!\!{_H}} \mathrm{d}r\; 
= \, \biggl( \displaystyle \int_{0}^{t} G^{n}\bigl(\overline{v}_{n}(r)\bigr) \, \mathrm{d}\overline{W}^{n}_r \, , z \biggr)_{\!\!\!{_H}} 
\end{equation}
\[
J_1^{n} + J_{2}^{n} + J_{3}^{n} + J_{4}^{n} = J_{5}^{n} \; .
\]
We pass then to the limit in \eqref{PbVar1} for almost all $t\in [0,T]$.

\noindent
For $J^{1}_n$, on the one hand we estimate  
\[
|  \bigl( \overline{v}_n(t) \, , \, z \bigr)_{\!\!{_H}} - \left( \overline{v}(t) \, , \, z \right)_{_H} | 
\leq \| \overline{v}_n - \overline{v} \|_{C^{0}([0,T] ; \D(A^{-3/2})\,)} \, \| z \|_{\D(A^{3/2})}\; . 
\]
On the other hand, since the restriction of $P_n$ to $H$ is the orthogonal projection onto $H_n$, we have
\begin{align*}
| & \bigl( P_n \overline{v}_n(0) \, , \, z \bigr)_{\!{_H}} \, - \, \bigl( \overline{v}(0) \, , \, z \bigr)_{\!{_H}} | \\
&= \bigl |  \bigl( \overline{v}_n(0) \, , \, P_n(z) - z \bigr)_{\!\!{_H}} - \bigl( \overline{v}(0) - \overline{v}_{n}(0) \, , \, z \bigr)_{\!\!{_H}} \bigr | \\
& \leq \bigl | \overline{v}_n(0) \bigr|_{_H} \,  \left| P_n(z) - z \right|_{_H} + \left\| \overline{v} - \overline{v}_{n} \right\|_{_{C^{0}([0,T] ; \D(A^{-3/2})\,)}} \, \left\| z \right\|_{_{\D(A^{3/2})}}.
\end{align*}
By \eqref{Estimevn1}, we get 
\[
\overline{\mathbb{E}}\left[ \, \left| \overline{v}_n(0) \right|_{_H} \,  \left| P_n(z) - z \right|_{_H} \, \right] \leq C \, \, \left| P_n(z) - z \right|_{_H} \; .
\]
We apply \eqref{ConvergenceProjection} and  \eqref{Convergence} and we obtain by thinning the sequence that 
\[
J_1^{n} \underset{n\to+\infty}{\longrightarrow} \bigl(\overline{v}(t) - \overline{v}(0)\, , \, z\bigr), \quad  \overline{\mathbb{P}}\,  \text{a.s} \; .
\]
\noindent
For $J_{2}^{n}$, since $A$ is symmetric for the inner product of $H$, we have 
\[
 \left| \int_{0}^{t} \left( A\overline{v}_{n}(r) \, , \, z \right)_{_H} \mathrm{d}r  - \int_{0}^{t} \left( A\overline{v}(r) \, , \, z \right)_{_H} \mathrm{d}r \right| \leq  \left\| \overline{v}_n - \overline{v}\right\|_{_{L^{2}([0,T] ;H)}} \; \left| A\, z \right|_{_H} 
\]
in particular, by (\ref{Convergence}), we obtain $
J_{2}^{n} \underset{n\to+\infty}{\longrightarrow} \displaystyle \int_{0}^{t} \bigl( A\overline{v}(r) \, , \, z \bigr)_{\!{_H}} \mathrm{d}r $,  \quad  $\overline{\mathbb{P}} \, \text{a.s}$ . \\
\noindent
For $J_{3}^{n}$, in the same way as for $J_{1}^{n}$, we have,
\begin{align*}
& \left| \int_{0}^{t} \bigl(B^{n} \overline{v}_{n}(r) \, , \, z \bigr)_{\!{_H}} \mathrm{d}r  - \int_{0}^{t} \bigl(B \overline{v}(r) \, , \, z \bigr)_{\!{_H}} \mathrm{d}r \right| \\
&\leq  \left| \int_{0}^{t} \bigl( B \overline{v}_{n}(r) \, , \, P_n(z) - z \bigr)_{\!{_H}} \mathrm{d}r \, \right| \, + \, \biggl| \int_{0}^{t} \bigl( B \overline{v}_{n}(r) - B \overline{v}(r) \, , \, z \bigr)_{\!{_H}} \mathrm{d}r \, \biggr| \\
&:= \, J_{3 \, , \, 1}^{n} \, + \, J_{3 \, , \, 2}^{n} \; .
\end{align*}
For $J_{3 \, , \, 1}^{n}$, we have $P_n(z) -z \in V$ since $\D(A^{3/2}) \subset V$. We infer by \eqref{BContinuite1}, 
\begin{align*}
 \overline{\mathbb{E}}\left[ \, {|}  \int_{0}^{t} \left( B \overline{v}_{n}(r) \, , \, P_n(z) - z \right)_{_H} \, \mathrm{d}r  \, {|}  \, \right]
& \leq \,  \overline{\mathbb{E}}\left[ \, \int_{0}^{t} \|  B \overline{v}_{n}(r) \|_{_{V'}} \; \|  P_nz - z  \|_{_{V}}  \, \mathrm{d}r \, \right] \\
&\leq \,  C \, \overline{\mathbb{E}}\left[ \, \int_{0}^{t}  \| \overline{v}_{n}(r) \|^{2}_{_{V}}\, \mathrm{d}r \, \right] \; \|  P_nz - z  \|_{_{V}}\; .
\end{align*}
By \eqref{Estimevn1}, we have by thinning the sequence, $J_{3 \, , \, 1}^{n} \underset{n\to+\infty}{\longrightarrow} 0$, \; $\overline{\mathbb{P}} \, \text{a.s}$. For $J_{3 \, , \, 2}^{n}$, we have
\[
\int_{0}^{t} \bigl( B \overline{v}_{n}(r) - B \overline{v}(r) \, , \, z \bigr)_{\!{_H}} \mathrm{d}r \,  =  \int_{0}^{t} \bigl\langle  B \overline{v}_{n}(r) - B \overline{v}(r) \, , \, z \bigr\rangle_{\D(A^{-\gamma})\times \D(A^{\gamma})} \mathrm{d}r  .
\]
The bilinearity of $B$ implies for all $r\in [0,T]$, 
\[
B \overline{v}_{n}(r) - B \overline{v}(r) = B \bigl(\overline{v}_{n}(r) - \overline{v}(r) \, , \, \overline{v_{n}}(r)  \bigr) \; {+} \;  B \bigl( \overline{v}(r) \, , \, \overline{v}_{n}(r) - \overline{v}(r)  \bigr) \; . 
\]
Due to \eqref{EstimB} and the Cauchy-Schwarz inequality, we obtain also
\begin{align*}
& \left| \int_{0}^{t} \left( B \overline{v}_{n}(r) - B \overline{v}(r) \, , \, z \right)_{_H} \mathrm{d}r \, \right| \\
&  \leq C \int_{0}^{t} | \overline{v}_n(r) - \overline{v}(r) |_{_H} \, | \overline{v}_n(r) |_{_H} \, \| z \|_{_{\D(A^{\gamma})}} \, \mathrm{d}r \; + \; C \int_{0}^{t} | \overline{v}(r) |_{_H} \, | \overline{v}_n(r) - \overline{v}(r) |_{_H} \, \| z \|_{_{\D(A^{\gamma})}} \, \mathrm{d}r \\
&  \leq C\, \| z \|_{_{\D(A^{\gamma})}} \, \left[ \sup_{n} \| \overline{v}_{n} \|_{_{L^{2}([0,T] ; H) }} \, + \, \| \overline{v} \|_{_{L^{2}([0,T] ; H) }}  \right] \; \| \overline{v}_n - \overline{v} \|_{_{L^{2}([0,T] ; H) }} .
\end{align*}
 As $(\overline{v}_{n})_n$ converges to $\overline{v}$ in $L^{2}([0,T] ; H), \; \overline{\mathbb{P}} \;  \text{a.s}$, we conclude by thinning the sequence that $J_{3 \, , \, 2}^{n} \underset{n\to+\infty}{\longrightarrow} 0$, \; $\overline{\mathbb{P}} \, \text{a.s}$. We get hence 
\[
J_{3}^{n} \underset{n\to+\infty}{\longrightarrow} \displaystyle \int_{0}^{t} \left( B \overline{v}(r) \, , \, z \right)_{_H} \mathrm{d}r , \qquad  \overline{\mathbb{P}} \, \text{a.s} \;. 
\]
For $J_{4}^{n}$, it can be observed that for almost all  $t\in [0,T]$,
\begin{align*}
&\left| \int_{0}^{t} \bigl( F^{n} \overline{v}_{n}(r) \, , \, z \bigr)_{\!{_H}} \mathrm{d}s  - \int_{0}^{t} \bigl( F \overline{v}(r) \, , \, z \bigr)_{_H} \mathrm{d}r \right| \\
&\leq \biggl| \int_{0}^{t} \bigl( F \overline{v}_{n}(r) \, , \, P_n(z) - z  \bigr)_{\!{_H}} \mathrm{d}r \biggr|\, +  \, \biggl|\int_{0}^{t} \bigl( F \overline{v}_{n}(r) - F \overline{v}(r) \, , \, z  \bigr)_{_H} \mathrm{d}r \biggr|  \\
&:= J_{4,1}^{n} \; + \; J_{4,2}^{n}.
\end{align*}
For $J_{4,1}^{n}$, as $P_n(z) -z \in V$, we have by \eqref{EstimF1}, 
\begin{align*}
\left| J_{4,1}^{n} \right|
&\leq \int_{0}^{t}  \| F\bigl(\overline{v}_{n}(r)\bigr) \|_{_{V'}} \, \mathrm{d}r \quad  \| P_nz - z \|_{_{V}} \\
&\leq C\, \varepsilon^{2} \,\bigl(\| \overline{v}_{n} \|_{_{L^{2}([0,T] ; V)}} + 1 \bigr) \; \| P_nz - z \|_{_{V}} \; .
\end{align*}
By \eqref{Estimevn1}, we obtain by thinning the sequence that $J_{4 \, , \, 1}^{n} \underset{n\to+\infty}{\longrightarrow} 0$, \; $\overline{\mathbb{P}}$ a.s . Furthermore, for $J_{4,2}^{n}$, the continuity of $F$ \eqref{EstimF2}, implies that
\begin{align*}
\left| J_{4,2}^{n} \right|
&\leq \int_{0}^{T}  \| F\bigl(\overline{v}_{n}(r)\bigr) - F\bigl(\overline{v}(r)\bigr) \|_{_{\D(A^{-1})}} \, \mathrm{d}r \quad  \| z \|_{_{\D(A)}} \\
&\leq C\, \varepsilon^{2} \,  \| \overline{v}_{n} - \overline{v} \|_{_{L^{2}([0,T] ; H)}} \; \| z \|_{_{\D(A)}} \; .
\end{align*}
Due to (\ref{Convergence}), by thinning the sequence we infer that $J_{4 \, , \, 2}^{n} \underset{n\to+\infty}{\longrightarrow} 0$, \; $\overline{\mathbb{P}}$ a.s. We obtain also that 
\[
J_{4}^{n} \underset{n\to+\infty}{\longrightarrow} \displaystyle \int_{0}^{t} \bigl( F \overline{v}(r) \, , \, z \bigr)_{\!{_H}} \mathrm{d}r,  \qquad  \overline{\mathbb{P}} \, \text{a.s}. 
\]
We shall study the convergence of the martingale term $J_{5}^{n}$ using lemma \ref{ConvSto}. For a.s $(t,\omega)\in [0,T] \times \overline{\Omega}$, we have,
\begin{align}
&\| G^{n}\bigl(\overline{v}_{n}(t)\bigr) - G\bigl(\overline{v}(t)\bigr) \|^{2}_{\mathcal{L}_{2}(L^2(\mathcal S,{\mathbb R}^d) , V')} \nonumber
\\
&\leq 2 \left( \| G^{n}\bigl(\overline{v}_{n}(t)\bigr) - G^{n}\bigl(\overline{v}(t)\bigr) \|^{2}_{\mathcal{L}_{2}(L^2(\mathcal S,{\mathbb R}^d) , V')} + \| G^{n}(\overline{v}(t)) - G(\overline{v}(t)) \|^{2}_{\mathcal{L}_{2}(L^2(\mathcal S,{\mathbb R}^d) , V')} \right) \label{Ineg1} 
\end{align}
On the one hand, due to \eqref{ContinuiteG2}, we get
\[
\| G^{n}(\overline{v}_{n}) - G^{n}(\overline{v}) \|^{2}_{_{L^{2}([0,T] ; \mathcal{L}_{2}(L^2(\mathcal S,{\mathbb R}^d) , V') \, )}}  \leq C \varepsilon^{2} \, \| \overline{v}_{n} - \overline{v} \|^{2}_{_{L^{2}([0,T] ; H)}}
\]
in particular we infer with (\ref{Convergence}) the convergence in probability
\begin{equation}\label{Etude1}
G^{n}(\overline{v}_{n}) - G^{n}(\overline{v}) \underset{n\to+\infty}{\longrightarrow} 0 \quad \text{in} \quad  L^{2}\bigl([0,T] \, ; \, \mathcal{L}_{2}(L^2(\mathcal S,{\mathbb R}^d) , V') \, \bigr) \; . 
\end{equation}
On the other hand, by  (\ref{ProjProp1}), it can be observed 
\begin{align*}
\| G^{n}(\overline{v}(t)) - G(\overline{v}(t)) \|^{2}_{{\mathcal{L}_{2}(L^2(\mathcal S,{\mathbb R}^d) , V') \, )}}  
&\leq \| (P^{n} - I) G(\overline{v}(t)) \|^{2}_{{ \mathcal{L}_{2}(L^2(\mathcal S,{\mathbb R}^d) \,, \, H) \, )}}   \\
&\leq C \, \lambda_{n+1}^{{-1}} \; \| G(\overline{v}(t)) \|^{2}_{_{\mathcal{L}_{2}(L^2(\mathcal S,{\mathbb R}^d) , V')}}\\
&\leq C \, \varepsilon^{2} \; \lambda_{n+1}^{{-1}} \, \left( 1 + | \overline{v}(t) |^{2}_{_H} \, \right) \; ,
\end{align*}
in particular, we have also 
\[
\overline{\mathbb{E}}\left[ \, \left| \int_{0}^{T} \| G^{n}(\overline{v}(t)) - G(\overline{v}(t)) \|^{2}_{_{\mathcal{L}_{2}(L^2(\mathcal S,{\mathbb R}^d) , V')}}  \mathrm{d}t \, \right|^{2} \, \right] \leq C \, \varepsilon^{2} \; \lambda_{n+1}^{{-1}} \, \left( 1 + \overline{\mathbb{E}} \left[ \sup_{0\leq t \leq T} | \overline{v}(t) |^{4}_{_H} \, \right] \, \right)
\]
which converges to $0$ when $[n\to\infty]$ since \eqref{EStimV1} is valid and $\lambda_{n}^{{-1}} \to 0$. We obtain in particular the following convergence in probability
\begin{equation}\label{Etude2}
G^{n}(\overline{v}) - G(\overline{v}) \underset{n\to+\infty}{\longrightarrow} 0 \quad \text{in} \quad  L^{2}\bigl([0,T] \, ; \, \mathcal{L}_{2}(L^2(\mathcal S,{\mathbb R}^d) , V') \, \bigr) \; . 
\end{equation}
Combining \eqref{Ineg1}, \eqref{Etude1} and \eqref{Etude2}, we apply Lemma \ref{ConvSto} to  obtain the convergence in probability
\begin{equation}\label{Conv3}
  \int_{0}^{t} G^{n}\bigl(\overline{v}_{n}(r)\bigr) \, \mathrm{d}\overline{W}_{r}^{n} \underset{n\to+\infty}{\longrightarrow} \int_{0}^{t} G\bigl(\overline{v}(r)\bigr) \, \mathrm{d}\overline{W}_{r} \quad \text{in} \quad  L^{2}([0,T] \, ; \, V') \; .  
\end{equation}
Furthermore, we shall show that $\left( \int_{0}^{t} G^{n}(\overline{v}_{n}(r)) \, \mathrm{d}\overline{W}^{n}_r  \right)_n$ is uniformly integrable: in fact, we apply the BDG inequality (proposition \ref{BDG})
\[
 \overline{\mathbb{E}} \left[ \, \left| \int_{0}^{T} \left\| \int_{0}^{t} G^{n}(\overline{v}_{n}(r)) \, \mathrm{d}\overline{W}^{n}_r  \right\|^{2}_{V'} \, \mathrm{d}t \, \right| \, \right] 
\leq C  \, \overline{\mathbb{E}} \left[ \, \int_{0}^{T} \left\| G^{n}(\overline{v}_{n}(t)) \, \right\|^{2}_{\mathcal{L}_{2}(L^2(\mathcal S,{\mathbb R}^d) , V')} \, \mathrm{d}t \, \right].
\]
By \eqref{ContinuiteG2} and \eqref{Estimevn1}, we have
\begin{align*}
 \overline{\mathbb{E}} \left[ \, \left| \int_{0}^{T} \left\| \int_{0}^{t} G^{n}(\overline{v}_{n}(r)) \, \mathrm{d}\overline{W}^{n}_r  \right\|^{2}_{V'} \, \mathrm{d}t \, \right| \, \right]  
&\leq C\, \varepsilon^{2} \; \overline{\mathbb{E}} \left[ \, \int_{0}^{T}  ( 1 + \left| \overline{v}_{n}(t) \right|^{2}_{_H}) \mathrm{d}t \, \right] \\
&\leq C\, \varepsilon^{2}\; .
\end{align*}
By the Vitali convergence theorem, the convergence (\ref{Conv3}) is also true in $L^{1}\bigl(\overline{\Omega} , L^{2}([0,T] ; V')\,\bigr)$. By thinning the sequence, we have the convergence  
\[
J_{5}^{n} \underset{n\to+\infty}{\longrightarrow} \left( \displaystyle \int_{0}^{t} G(\overline{v}(r)) \, \mathrm{d}\overline{W}_r \; , \, z \right)_{_H} \qquad \overline{\mathbb{P}} \, \text{a.s} .
\]
\noindent
Thus, we pass to the limit $[n\to\infty]$ in (\ref{PbVar1}) to obtain for all $z \in \D(A^{3/2})$ and almost surely $(t,\omega)\in [0,T]\times\overline{\Omega}$ that
\begin{equation*}
\bigl( \overline{v}(t) - v_0 \, , \, z\bigr)_{\!{_H}} \, + \, \displaystyle \int_{0}^{t} \left( A\overline{v}(r) \, , \, z \right)_{\!{_H} }\mathrm{d}r \, + \, \displaystyle \int_{0}^{t} \bigl( B \overline{v}(r) \, , \, z \bigr)_{\!{_H}} \, \mathrm{d}r  
\end{equation*}
\begin{equation*}
+ \, \displaystyle \int_{0}^{t} \bigl( F \overline{v}(r) \, , \, z\bigr)_{\!{_H}} \mathrm{d}r \; 
= \, \left( \displaystyle \int_{0}^{t} G(\overline{v}(r)) \, \mathrm{d}\overline{W}_r \, , \, z \right)_{\!\!{_H}} \; . 
\end{equation*}
 By density of $\D(A^{3/2})$ in $H$, $\overline{v}$ verify for almost surely $(t,\omega)\in [0,T] \times \overline{\Omega}$
\begin{equation}\label{EqSolution}
\overline{v}(t) - v_0  +  \displaystyle \int_{0}^{t} \bigl(A\overline{v}(r) +B \overline{v}(r) +F \overline{v}(r)\bigr) \, \mathrm{d}r
=\displaystyle \int_{0}^{t} G\bigl(\overline{v}(r)\bigr) \, \mathrm{d}\overline{W}_r 
\end{equation}
in $H$. We recall that $\overline{v} \in  L^{2}\bigl(\overline{\Omega} \, ; \, L^{2}([0,T] , V)\bigr) \cap L^{2}\bigl(\overline{\Omega} \, ; \, L^{\infty}([0,T] , H)\bigr)$. {It follows} that $\overline{v}$ is a martingale solution of \eqref{AbstractProblem}.

\subsection{The case of $2D$ fluid flow}\label{Case2D}
When $d=2$, we prove in this section that the solutions of the equation \eqref{EqSolution} are {in $C([0,T],H)$ and are }unique. 

{We first prove that $\overline{v} \in C^{0}([0,T] \, , \, H)$ \, $\overline{\mathbb{P}} $ a.s. To that end, we consider the following equation with $\overline{v}$ ,
\begin{equation}\label{EquationInterm}
\left\{
    \begin{array}{ll}
        d_t z(t) \; + \; Az(t)  \, \mathrm{d}t \; = \; G\bigl(\overline{v}(t)\bigr)\, \mathrm{d}\overline{W}_t,  &  \\
        z(0)=\overline{v}_0, & \;
    \end{array}
\right.
\end{equation}
which admits a unique solution $z\in L^{2}\bigl(\overline{\Omega} \, ; \, C^{0}([0,T] , H)\bigr)$. Furthermore, by the Ito formula applied to $x\to \frac{1}{\,2\,} | x|^{2}_{H} $, the solution $z$ of \eqref{EquationInterm} verifies 
\[
\overline{\mathbb{E}}\left[ \, \sup_{0\leq r \leq T} \,  |z(r)  |^{2}_{_H} \; + \; \displaystyle \int_{0}^{T} \| z(r) \|^{2}_{_V}  \, \mathrm{d}r \, \right] < \infty \;, 
\]
in particular $z \in L^{2}([0,T] \, ; \, V)$ \, $\overline{\mathbb{P}}$ a.s. Let $\Tilde{u}(t) := \overline{v}(t) - z(t)$ then $\Tilde{u}$ belongs to $L^{2}([0,T] \, ; \, V)$ \, $\overline{\mathbb{P}}$ a.s. and is solution of
\begin{equation*}
\left\{
    \begin{array}{ll}
        d_t \Tilde{u}(t) \; + \; A\Tilde{u}(t)  \, \mathrm{d}t \; + B\bigl(\Tilde{u}(t)+z(t))\, \mathrm{d}t \; + \; F(\Tilde{u}(t)+z(t)\bigr)\, \mathrm{d}t \; = \; 0  &  \\
        \Tilde{u}(0)=v_0 & 
    \end{array}
\right.
\end{equation*}
Using the hypothesis on $A$, $B$ and $F$ and the energy estimates, $d_t \Tilde{u}$ belongs to $L^{2}([0,T] \, ; \, V')$ \, $\mathbb{P}$ a.s., we conclude that $\Tilde{u}\in C([0,T] ; H)$ \, $\overline{\mathbb{P}}$ a.s. \\
\noindent
}

Concerning uniqueness, let us consider $v_1$ and $v_2$ two solutions of \eqref{AbstractProblem} on the same probability space $(\Omega, \mathcal{F} , (\mathcal{F}_t)_t , \mathbb{P} )$. We shall prove that $\mathbb{P}$ a.s, \; $v_1(t)=v_2(t)$ for all $t\in [0,T]$.
\\ \ \\
Let $\Tilde{V}(t) := v_1(t) -v_2(t)$ for $t\in[0,T]$, $\Tilde{V}$ satisfies the equation
\begin{multline*}
d_t \Tilde{V}(t) \; + \; \left(  A\Tilde{V}(t) + B(v_1)(t) - B(v_2)(t) +  F(v_1)(t) - F(v_2)(t)  \right) \mathrm{d}t \\
= \; \bigl( G(v_1)(t) - G(v_2)(t) \bigr) \, \mathrm{d}W_t  \; . \qquad 
\end{multline*}
Let $g(t):= \alpha \int_{0}^{t} \| v_2(r) \|^{2}_{_V} \, \mathrm{d}r$ and $e(t) := \exp\left(-g(t)\,\right)$ with $\alpha$ a positive constant chosen below. We apply the Ito formula\footnote{The use of Ito formula is not fully rigorous here. A regularization argument should be used. For instance, we can write the equation satisfied by $P_n\tilde V$, apply Ito formula to $e(t)|P_n\tilde V|^2$ and then let $n\to\infty$.} to $(t,x) \to \frac{1}{\,2\,} e(t) \, | x |^{2}_{_H}$ for $t\in [0,T]$ and $x\in H$ and obtain for $t\in [0,T]$, 
\begin{multline}
\frac{1}{\,2\,} e(t) \, | \Tilde{V}(t) |^{2}_{_H}  = 
\displaystyle \int_{0}^{t} e(r) \left( \,\Tilde{V}(r) \; , \; \left[ G(v_1(r)) - G(v_2(r)) \right]\, \mathrm{d}W_r \right)_{_H}\\
-\; \dfrac{1}{\,2\,} \; \displaystyle \int_{0}^{t} g'(r) \; e(r) \; | \Tilde{V}(r) |^{2}_{_H} \, \mathrm{d}r 
\\
- \displaystyle \int_{0}^{t} e(r) \, \left(  A\Tilde{V}(r) + B\bigl(v_1(r)\bigr) - B\bigl(v_2(r)\bigr) +  F\bigl(v_1(r)\bigr) - F\bigl(v_2(r)\bigr) \, , \, \Tilde{V}(r) \right)_{\!\!{_H}} \, \mathrm{d}r \; 
\\
+ \frac{1}{\,2\,} \displaystyle \int_{0}^{t} e(r) \,  \| G\bigl(v_1(r)\bigr) - G\bigl(v_2(r)\bigr)\|^{2}_{\mathcal{L}_{2}( L^2(\mathcal S,{\mathbb R}^d) , H)} \, \mathrm{d}r.
\end{multline}
We have for all $r\in[0,T]$, by \eqref{RelationA} and bilinearity of $B$, that
\[
( A \Tilde{V}(r) \, , \, \Tilde{V}(r) )_{_H} = \dfrac{1}{\,R_e\,} \, \| \Tilde{V}(r) \|^{2}_{_{V}}
\]
\[
B\bigl(v_1(r)\bigr) - B\bigl(v_2(r)\bigr) = B\bigl(v_1(r) , \Tilde{V}(r)\bigr) + B\bigl(\Tilde{V}(r) , v_2(r)\bigr).
\]
In particular by \eqref{PropB1}, we obtain
\begin{align*}
\left(  B(v_1(r)) - B(v_2(r)) \; , \; \Tilde{V}(r) \right)_{_H} 
&=  \left( \, B(\Tilde{V}(r) , v_2(s) )\; \; , \; \Tilde{V}(r) \right)_{_H} \\
&= \, b\bigl(\Tilde{V}(r) , v_2(s) , \Tilde{V}(r) \bigr) \; .
\end{align*}
We apply the Cauchy-Schwarz inequality and we get
\[
\bigl | b\bigl(\Tilde{V}(r) , v_2(s) , \Tilde{V}(r) \bigr)\bigr | \; \leq \; \| \Tilde{V}(r) \|^{2}_{_{L^{4}(\mathcal{S})}} \, \| v_{2}(r) \|_{_{V}} \; .
\]
By Gagliardo-Nirenberg inequality we deduce 
\[
\biggl |  \biggl(  B\bigl(v_1(r)\bigr) - B\bigl(v_2(r)\bigr) \; , \, \Tilde{V}(r)  \biggr)_{\!\!{_H}} \biggr| 
\leq C  \, | \Tilde{V}(r) |_{\!\!{_H}} \, \|v_{2}(r) \|_{_{V}} \, \| \Tilde{V}(r) \|_{_{V}}
\]
and we obtain thus,  
\begin{equation}\label{UniciteInterm}
\biggl |  \biggl(  B\bigl(v_1(r)\bigr) - B\bigl(v_2(r)\bigr) \; , \, \Tilde{V}(r)  \biggr)_{\!\!{_H}} \biggr| 
\leq \dfrac{1}{\,2 R_e \,} \, \| \Tilde{V}(r) \|^{2}_{_{V}} \; + \; C_{1} \, | \Tilde{V}(r) |^{2}_{\!{_H}}\|v_{2}(r) \|^{2}_{_{V}} \; .
\end{equation}
Furthermore, by definition of $F$ (and $a$), we have 
\begin{align}
\biggl(  \Tilde{V}(r) \, ,  &  \,  F\bigl(v_1(r)\bigr) - F\bigl(v_2(r)\bigr) \,  \biggr)_{\!\!{_H}} \nonumber\\ 
&= \, \varepsilon^{2} \left( \Tilde{V}(r) \; , \;  \Tilde{V}(r)\bcdot \nabla u_s \, \right)_{\!\!{_H}}  - \dfrac{\, \varepsilon^{2}\; }{2} \, \bigl( \Tilde{V}(r) \, , \,  \nabla \bcdot [\, a\nabla \Tilde{V}(r)\, ] \, \bigr)_{\!{_H}} \nonumber \\
&= \, \varepsilon^{2} \left( \Tilde{V}(r) \; , \;  \Tilde{V}(r)\bcdot \nabla u_s \, \right)_{\!{_H}} \; + \; \dfrac{\, \varepsilon^{2}\; }{2} \, \displaystyle \sum_{k=0}^{\infty} \left| (\phi_{k} \bcdot \nabla) \Tilde{V}(r) \right|^{2}_{_H} \nonumber \\
&\leq  \, C \, \varepsilon^{2} | \Tilde{V}(r) |^{2}_{\!{_H}} \; + \; \dfrac{\, \varepsilon^{2}\; }{2} \, \displaystyle \sum_{k=0}^{\infty} \left| (\phi_{k} \bcdot \nabla) \Tilde{V}(r) \right|^{2}_{{_H}} \label{IntermEtudeFUnicite}
\end{align}
since \eqref{hyp} is satisfied. The second part of the last term is exactly equal to 
\[
\dfrac{1}{\,2\,} \,  \| G(v_1(r)) - G(v_2(r))\|^{2}_{\mathcal{L}_{2}( L^2(\mathcal S,{\mathbb R}^d) , H)}  .
\]
We choose $\alpha:= 2 C_1$ with $C_1$ defined in \eqref{UniciteInterm} and obtain
\begin{align}
\displaystyle \dfrac{1}{\,2\,}  e(t)  | \Tilde{V}(t) |^{2}_{_H} \;  + \; & \dfrac{1}{\, 2 R_e \, }  \int_{0}^{t} e(r) \| \Tilde{V}(r) \|^{2}_{_V}  \, \mathrm{d}r 
\leq \displaystyle C \, \varepsilon^{2} \,  \int_{0}^{t}  e(r)  | \Tilde{V}(r) |^{2}_{_H} \mathrm{d}r \nonumber \\
&+ \;  \displaystyle \int_{0}^{t}  e(r) \left(  \Tilde{V}(r) \, , \, \left[ G(v_1(r)) - G(v_2(r)) \right] \, \mathrm{d}W_r \right)_{\!{_H}}\; . \label{UniciteInterm0} 
\end{align}
As $\left( \displaystyle \int_{0}^{t} e(r) \left( \,\Tilde{V}(r) \; , \; \left[ G\bigl(v_1(r)) - G(v_2(r)\bigr) \right]\, \mathrm{d}W_r \right)_{\!{_H}} \right)_{t}$
is a martingale. We apply  Gronwall's lemma to obtain for all $t\in [0,T]$ 
\begin{equation*}
\mathbb{E}[ \, e(t) \, | \Tilde{V}(t) |^{2}_{_H} \, ]= 0 . 
\end{equation*}
We deduce that for all $t\in [0,T]$, $\Tilde V(t)=0$ $\mathbb{P}$ a.s. Since $\Tilde V$ is continous, we deduce that, $\mathbb{P}$ a.s, $\Tilde V(t)=0$  for all 
$t\in [0,T]$. This proves uniqueness. 

Using an argument due to Gyongy and Krylov (see for instance \cite{Debussche-et-al-2011}, section 5), we conclude that the whole sequence $(v_n)_n$ converges to a unique solution of \eqref{EqSolution} in probability in the original stochastic basis. This in particular gives the existence of a probabilistic strong solution in the $2D$ case.

\section{The limit $\varepsilon \to 0$}\label{PassageLimiteEpsilon}
Let $v_0\in H$. For all $\varepsilon >0$, we have proved that the abstract problem \eqref{EqNVS3} admits martingale solutions: we have built a family $(v_\varepsilon)_{\varepsilon>0}$ of solutions. We show now that $(v_{\varepsilon})_{\varepsilon>0}$ converges when $[\varepsilon\to 0^{+}]$ to a solution $v$ of the following deterministic Navier-Stokes equation (see section \ref{NVSdeterministe}) 
\begin{equation}\label{NVSdeterministEQ}
\left\{
    \begin{array}{ll}
        d_t v(t) \; + \; \dfrac{1}{\, R_e\, } \Delta v(t)  \, \mathrm{d}t \; + \; \bigl(v(t) \bcdot \nabla \bigr) v(t)\, \mathrm{d}t \;  = \; 0  &  \\
        \nabla \bcdot v = 0 &  \\
        v(0)=v_0 & 
    \end{array}
\right.
\end{equation}
with $v_0 \in H$. This equation can be written as the following abstract problem
\begin{equation}\label{NVSdeterEQabstrac}
\left\{
    \begin{array}{ll}
        d_t v(t) \; + \; A v(t)  \, \mathrm{d}t \; + \; B v(t) \, \mathrm{d}t \;  = \; 0   \\
        v(0)=v_0 \; .
    \end{array}
\right.
\end{equation}

\subsection{$2D$ case}\label{Convergence2D}
In $2D$ we have seen that or all $\varepsilon>0$ the stochastic equation \eqref{AbstractProblem} admits a unique  solution $v_{\varepsilon}$ associated to the stochastic basis $(\Omega , \mathcal{F}, (\mathcal{F}_t)_{t\in[0,T]}, \mathbb{P}  )$ and the cylindrical Wiener process $W$. The deterministic equation \eqref{NVSdeterEQabstrac} admits an unique weak solution $v$. We prove below that $(v_{\varepsilon})_{\varepsilon}$ converges to $v$ in the following sense 
\[ 
\mathbb{E} \left[ \sup_{0\leq t \leq T}  e(t)\, |v_{\varepsilon}(t) - v(t) |^{2}_{_H}  \, \right] + \mathbb{E}\left[ \int_{0}^{T} e(r) \| v_{\varepsilon}(r) - v(r)\|^{2}_{_V}\, \mathrm{d}r \right]\underset{\varepsilon\to 0^{+}}{\longrightarrow} 0,
\] 
where $e(t) := \exp\left(- \alpha \int_{0}^{t} \| v(r) \|^{2}_{_V} \, \mathrm{d}r \,\right)$ and $\alpha >0$ is a positive constant that will be further specified in the following. 
{Since $v\in L^2([0,T];V)$, we have:
$$
0<e(T)\le e(t),\quad t\in [0,T],
$$
and we deduce:
\[ 
\mathbb{E} \left[ \sup_{0\leq t \leq T}\, |v_{\varepsilon}(t) - v(t) |^{2}_{_H}  \, \right] + \mathbb{E}\left[ \int_{0}^{T} \| v_{\varepsilon}(r) - v(r)\|^{2}_{_V}\, \mathrm{d}r \right]\underset{\varepsilon\to 0^{+}}{\longrightarrow} 0,
\] 
it follows that $v_\varepsilon$ converges to $v$ in  $L^2(\Omega; L^\infty([0,T];H)\cap L^2([0,T];V))$.}

For all $\varepsilon>0$, let $z_{\varepsilon}(t):= v_{\varepsilon}(t) - v(t)$ for all $t\in [0,T]$. The random variable $z_{\varepsilon}$  satisfies the following equation 
\begin{equation}\label{EquationAbstactLimitEpsilon}
\left\{
    \begin{array}{ll}
        d_t z_{\varepsilon}(t) \; + \; Az_{\varepsilon}(t)  \, \mathrm{d}t  + \bigl(B v_{\varepsilon}(t)\, - B v(t)  \bigr) \mathrm{d}t  +  F_{\varepsilon} v_{\varepsilon}(t)\, \mathrm{d}t \; = \; G_{\varepsilon}\bigl( v_{\varepsilon}(t)\bigr)\, \mathrm{d}W_t  &  \\
        z_{\varepsilon}(0)=0  \; .& 
    \end{array}
\right.
\end{equation}
Let $h(t):= \alpha \int_{0}^{t} \| v(r) \|^{2}_{_V} \, \mathrm{d}r$ and $e(t) = \exp\left(-h(t)\right)$. We apply the Ito formula to the function
$F(t,x) = \frac{1}{\,2\,} e(t) \, | x |^{2}_{_H}$ for $t\in [0,T]$ and $x\in H$  and we obtain for all $t\in [0,T]$, 
\begin{multline}
\frac{1}{\,2\,} e(t) \, | z_{\varepsilon}(t) |^{2}_{_H}  = 
\displaystyle \int_{0}^{t} e(r) \, \biggl(   z_{\varepsilon}(r) \, , \, G_{\varepsilon}\bigl(v_{\varepsilon}(r)\bigr) \, \mathrm{d}W^{\varepsilon}_r \biggr)_{\!{_H}} \, - \, \displaystyle \int_{0}^{t} h'(r) \; e(r) \; | z_{\varepsilon}(r) |^{2}_{\!{_H}} \, \mathrm{d}r \\
- \displaystyle \int_{0}^{t} e(r) \, \bigl(  A z_{\varepsilon}(r) + B v_{\varepsilon}(r) - B v(r) +  F_{\varepsilon} v_{\varepsilon}(r) \; , \; z_{\varepsilon}(r)\bigr)_{\!{_H}} \, \mathrm{d}r \; 
\\
+ \frac{1}{\,2\,} \displaystyle \int_{0}^{t} e(r) \,  \| G_{\varepsilon}(v_{\varepsilon}(r))\|^{2}_{{\mathcal{L}_{2}\bigl(L^{2}(\mathcal{S} , \mathbb{R}^d) \, , H \bigr)}} \, \mathrm{d}r.
\end{multline}
With the same arguments as in section \ref{Case2D}, we have for all $r\in [0,T]$
\[
\left|  \bigl(  B v_{\varepsilon}(r) - B v(r) \; , \, z_{\varepsilon}(r)  \bigr)_{\!{_H}} \right| 
\leq \dfrac{1}{\,3 R_e \,} \, \| z_{\varepsilon}(r) \|^{2}_{_{V}} \; + \; C_{1} \, | z_{\varepsilon}(r) |^{2}_{_H}\|v(r) \|^{2}_{_{V}},
\]
with $C_1>0$ a constant independent of $\varepsilon$. Furthermore, we have 
\[
\bigl( F_{\varepsilon} v_{\varepsilon}(r) \; , \; z_{\varepsilon}(r) \bigr)_{\!{_H}} = \left( F_{\varepsilon} v_{\varepsilon}(r) \; , \; v_{\varepsilon}(r) \right)_{\!{_H}} \; - \; \left( F_{\varepsilon} v_{\varepsilon}(r) \; , \; v(r) \right)_{\!{_H} } . 
\]
Due to \eqref{EstimF1} and Holder inequality, we infer 
\begin{align*}
| \left( F_{\varepsilon} v_{\varepsilon}(r) \; , \; v(r) \right)_{_H} | 
&\leq C\,\varepsilon^{2} \; \left( \|v_{\varepsilon}(r)\|_{_V} + 1\right) \, \|v(r)\|_{_V} \\ 
&\leq C\,\varepsilon^{2} \; \left( \|z_{\varepsilon}(r)\|_{_V}  + \|v(r)\|_{_V} + 1\right) \, \|v(r)\|_{_V} \\
&\leq \dfrac{1}{\,3R_e\,} \; \|z_{\varepsilon}(r)\|^{2}_{_V} \; + \;   C  \varepsilon^{2} \|v(r)\|^{2}_{_V}  \; + \; C \varepsilon^{2}.
\end{align*}
With the same arguments as in section \ref{PreuveEstimees}, we obtain 
\begin{align}
\left( F_{\varepsilon} v_{\varepsilon}(r) \; , \; v_{\varepsilon}(r) \right)_{_H} 
&\leq \dfrac{\, \varepsilon^{2}\; }{2} \, \displaystyle \sum_{k=0}^{\infty} \left| (\phi_{k} \bcdot \nabla) v_{\varepsilon}(r) \right|^{2}_{_H} \; +
\; C \, \varepsilon^{2}  \; |v_{\varepsilon}(t) |^{2}_{_H} + C\varepsilon^{2}\nonumber \\
&\leq \; \dfrac{\, \varepsilon^{2}\; }{2} \, \displaystyle \sum_{k=0}^{\infty} \left| (\phi_{k} \bcdot \nabla) v_{\varepsilon}(r) \right|^{2}_{_H} \; + \; C \, \varepsilon^{2}  \; |z_{\varepsilon}(r) |^{2}_{_H}  + C \, \varepsilon^{2}  \; |v(r) |^{2}_{_H} + C\varepsilon^{2}. \label{IntermFPassageLimite}
\end{align}
Finally, applying again the same arguments as in section \ref{PreuveEstimees} we obtain 
\begin{align*}
\dfrac{1}{\,2\,} \; & \| G_{\varepsilon}(v_{\varepsilon}(r))\|^{2}_{\mathcal{L}_{2}(L^2(\mathcal S,{\mathbb R}^d) , H))} \\
& \leq \frac{\,\varepsilon^{2}\,}{2}  \; \displaystyle \sum_{k=0}^{\infty} |(\phi_k \bcdot \nabla) v_{\varepsilon}(r)|^{2}_{_H} \; + \; C\varepsilon^{2} \; + \; 2 \varepsilon^{2}  |v_{\varepsilon}(r)|^{2}_{_H} \\
&\leq \frac{\,\varepsilon^{2}\,}{2}  \; \displaystyle \sum_{k=0}^{\infty} |(\phi_k \bcdot \nabla) v_{\varepsilon}(r)|^{2}_{_H} \; + \; C\varepsilon^{2} \; + \; 4 \varepsilon^{2}  |z_{\varepsilon}(r)|^{2}_{_H} \; + \; 4 \varepsilon^{2}  |v(r)|^{2}_{_H} 
\end{align*}
and the first term is exactly equal to a term of \eqref{IntermFPassageLimite}.\\
We choose $\alpha:= 2 C_1$. Since $|v(r)|^{2}_{_H} \leq C \|v(r)\|^{2}_{_V}$ for all $r\in [0,T]$. By the previous inequalities, we obtain for all $t\in[0,T]$ that
\begin{align}
&\displaystyle \dfrac{1}{\,2\,}  e(t)  | z_{\varepsilon}(t) |^{2}_{_H} \;  + \;  \dfrac{1}{\, 3 R_e \, }  \int_{0}^{t} e(r) \| z_{\varepsilon}(r) \|^{2}_{_V}  \, \mathrm{d}r \; \leq  \; C \, k_{\varepsilon}(t) \; + \; C \varepsilon^{2}\nonumber\\
&+\displaystyle C \varepsilon^{2}   \int_{0}^{t}  e(r)  | z_{\varepsilon}(r) |^{2}_{_H} \mathrm{d}r + \displaystyle \int_{0}^{t} e(r) \, \biggl(   z_{\varepsilon}(r) \, , \, G\bigl(v_{\varepsilon}(r)\bigr) \, \mathrm{d}W^{\varepsilon}_r \biggr)_{\!\!{_H}}, \label{UniciteInterm1} 
\end{align}
with the deterministic function $ \displaystyle k_{\varepsilon}(t) := \varepsilon^{2}   \int_{0}^{t}  e(r)  \| v(r) \|^{2}_{_V} \mathrm{d}r$. Taking expectation and using Gronwall lemma, we have also
\begin{equation}
    \mathbb{E}\left[ e(t) | z_{\varepsilon}(t)|^{2}_{_H} \right] \leq C  \left( k_{\varepsilon}(T) + \varepsilon^{2}\right)\; \exp(  C \varepsilon^{2} T),
\end{equation}
with $C>0$ a constant independent of $\varepsilon$. As $v\in L^{2}([0,T], V)$ we obtain 
\begin{equation}\label{InterAA}
\mathbb{E}\left[ e(t) | z_{\varepsilon}(t)|^{2}_{_H} \right] \underset{\varepsilon\to 0}{\longrightarrow} 0 \; \, \text{and}  \; \, \mathbb{E}\left[ \int_{0}^{t} e(r) \| z_{\varepsilon}(r)\|^{2}_{_V}\, \mathrm{d}r \right] \underset{\varepsilon\to 0}{\longrightarrow} 0 .
\end{equation}
Arguing as in section \ref{PreuveEstimees}, by the BDG inequality (proposition \ref{BDG}) and \eqref{UniciteInterm1} the following inequality holds
\begin{multline*}
\dfrac{1}{4} \, \mathbb{E}\left[ \sup_{t\in [0,T]} e(t) | z_{\varepsilon}(t)|^{2}_{_H} \right] 
\leq C\varepsilon^{2}  \, \mathbb{E}\left[ \displaystyle \int_{0}^{T} e(r) | z_{\varepsilon}(r)|^{2}_{_H} \, \mathrm{d}r\right] 
+ \; C \left( \varepsilon^{2}  +\varepsilon^{2} \int_{0}^{T} e(r)  \| v(r) \|^{2}_{_V} \mathrm{d}r \right).
\end{multline*}
As $v \in L^{2}([0,T] , V)$, we apply \eqref{InterAA} (since $V\subset H$) to obtain 
\[
\displaystyle \mathbb{E}\left[ \sup_{t\in [0,T]} e(t) | z_{\varepsilon}(t)|^{2}_{_H} \right] \underset{\varepsilon\to 0}{\longrightarrow} 0 . 
\]
\subsection{$3D$ case}
Let $(\varepsilon_n)_n$ a sequence which converge to $0$. For all $n\in \mathbb{N}$, by Theorem \ref{Theorem}, equation (\ref{AbstractProblem}) admits a martingale solution $v_{\varepsilon_{_n}}$ associated to the stochastic basis $(\Omega_{n} ,\mathcal{F}_{n}, (\mathcal{F}^{n}_t)_{t\in[0,T]}, \mathbb{P}_{n} )$ which verifies the following estimates 
\begin{align}
    \sup_{n\in\mathbb{N}} \, \mathbb{E}_{n} [ \, \sup_{r \in [0,T]} |v_{\varepsilon_n}(r)  |^{p}_{_H} \, ] &\leq C \qquad \forall \, p\geq 2, \label{Estimevepsilon1} \\
    \qquad \sup_{n\in\mathbb{N}} \, \mathbb{E}_{n} [ \,  \displaystyle \int_{0}^{T} \| v_{\varepsilon_n}(r) \|^{2}_{_V}  \, \mathrm{d}r \; ] &\leq C, \label{Estimevepsilon2}
\end{align}
by the same proof as \eqref{Estimevn1}. Arguing similarly to section \ref{Tight}, these estimates enable us to prove that 
$\left(\mathcal{L}(v_{\varepsilon_n}) \right)_{n}$ are tight in 
\[
L^{2}([0,T] \, ; \, H) \, \cap \, C^{0}([0,T] \, ; \, \D(A^{-3/2})\, ).
\]
By Skorohod embedding theorem, there exists a stochastic basis $(\overline{\Omega} ,\overline{\mathcal{F}}, (\overline{\mathcal{F}}_t)_{t\in[0,T]}, \overline{\mathbb{P}} )$ and $L^{2}([0,T] ; H) \, \cap \, C^{0}\bigl([0,T] ; \D(A^{-3/2})\, \bigr)$ valued random variables $(\overline{v}_{\varepsilon_n}, \overline{W}^{\varepsilon_{n}})$ and $(\overline{v} , \overline{W})$ such that  : $\overline{v}_{\varepsilon_n}$ has the same law of $v_{\varepsilon_n}$ on $L^{2}([0,T] ; H) \, \cap \, C^{0}\bigl([0,T] ; \D(A^{-3/2})\, \bigr)$ ;  $\overline{W}^{\varepsilon_{n}} \; n\geq 1$ and $\overline{W}$ are cylindrical Wiener processes with 
\begin{align}
\overline{v}_{\varepsilon_n} \to \overline{v} \quad  \text{in} \;  L^{2}([0,T] ; H) &\cap C^{0}\bigl([0,T] ; \D(A^{-3/2})\bigr) \qquad  \overline{\mathbb{P}} \;  \text{a.s}  \label{ConvergenceEpsilon} \\
\overline{W}^{\varepsilon_n} \to \overline{W} \quad  &\text{in} \; C^{0}([0,T] , U_{0}) \qquad \overline{\mathbb{P}} \;  \text{a.s} \label{ConvergenceEpsilonSto} \; .
 \end{align}
Each pair verify the following equation 
\begin{multline}\label{EqApproxEpsilon2}
\overline{v}_{\varepsilon_n}(t) - v_0 + \displaystyle \int_{0}^{t} \left[ A\overline{v}_{\varepsilon_n}(r) + B\overline{v}_{\varepsilon_n}(r) +F_{\varepsilon_n}\overline{v}_{\varepsilon_n}(r) \right] \mathrm{d}r = \displaystyle \int_{0}^{t} G_{\varepsilon_n}(\overline{v}_{\varepsilon_n}(r)) \mathrm{d}\overline{W}^{\varepsilon_n}_r.
\end{multline}
Arguing similarly to section \ref{PassageLimite}, it can be noticed that
\begin{align}
    \sup_{n\in\mathbb{N}} \, \overline{\mathbb{E}} [ \, \sup_{r \in [0,T]} |\overline{v}_{\varepsilon_n}(r)  |^{p}_{_H} \, ] &\leq C, \qquad \forall \, p\in \mathbb{N} ,\label{EstimevnEpsilon1} \\
    \qquad \sup_{n\in\mathbb{N}} \, \overline{\mathbb{E}} [ \,  \displaystyle \int_{0}^{T} \| \overline{v}_{\varepsilon_n}(r) \|^{2}_{_V}  \, \mathrm{d}r \; ] &\leq C \; . \label{EstimevnEpsilon2}
\end{align}
and $\overline{v}$ belongs to $L^{2}(\overline{\Omega} \, ; \, L^{2}([0,T] , V)) \cap L^{2}(\overline{\Omega} \, ; \, L^{\infty}([0,T] , H))$. We will pass to the limit $[n\to\infty]$ in \eqref{EqApproxEpsilon2} for almost $t\in [0,T]$. To that end, for all $n\in\mathbb{N} $, $t\in [0,T]$ and $z\in  \D(A^{3/2})$, we have
\begin{multline}
\label{eq-Kin}
\left( \overline{v}_{\varepsilon_n}(t) - v_0 \, , \, z\right)_{_H} \, + \, \displaystyle \int_{0}^{t} \left( A\overline{v}_{\varepsilon_n}(r) \, + \, B \overline{v}_{\varepsilon_n}(r) , \, z \right)_{_H} \mathrm{d}r \,  \\
+ \, \displaystyle \int_{0}^{t} \left( F_{\varepsilon_n} \overline{v}_{\varepsilon_n}(r) \, , \, z\right)_{_H} \mathrm{d}r \; 
= \, \biggl( \displaystyle \int_{0}^{t} G_{\varepsilon_n}\bigl(\overline{v}_{\varepsilon_n}(r)\bigr) \, \mathrm{d}\overline{W}^{\varepsilon_n}_r \, , z \biggr)_{\!\!{_H}}  
\end{multline}
\[
K_{1}^{\varepsilon_n} + K_{2}^{\varepsilon_n} + K_{3}^{\varepsilon_n}   = K_{4}^{\varepsilon_n} \; .
\]
As in the section \ref{PassageLimite}, for almost $t\in [0,T]$, by \eqref{EstimevnEpsilon1} and \eqref{EstimevnEpsilon2} it can be proved, by thinning the sequence, that 
\[
K_1^{\varepsilon_n} + K_{2}^{\varepsilon_n} \underset{n\to+\infty}{\longrightarrow} (\overline{v}(t) - v_{0}\, , \, z) \, +  \, \int_{0}^{t} \bigl( A\overline{v}(r) \, + \, B \overline{v}(r) \, , \, z \bigr)_{\!{_H}} \mathrm{d}r  \quad  \overline{\mathbb{P}} \, \text{a.s} .
\]  
Furthermore, by \eqref{EstimF1} we infer 
\[
\left| \int_{0}^{t} \bigl( F_{\varepsilon_n} \overline{v}_{\varepsilon_n}(r) \, , \, z\bigr)_{\!{_H}} \mathrm{d}r \right| \leq C\,\varepsilon_{n}^{2} \left( \int_{0}^{T} \| \overline{v}_{\varepsilon_n}(r)\|_{_V} \, \mathrm{d}r \; + \; 1 \, \right) \|z\|_{_V} .
\]
Due to \eqref{EstimevnEpsilon2} we obtain, by thinning the sequence, that 
\[
K_{3}^{\varepsilon_n} \underset{n\to+\infty}{\longrightarrow} 0   \quad  \overline{\mathbb{P}} \; \text{a.s} \, .
\]
Finally, by (\ref{EstimG0}) and by (\ref{EstimevnEpsilon2}) we have
\begin{align*}
\overline{\mathbb{E}}\left[ \, \left| \int_{0}^{T} \| G_{\varepsilon_{n}}(\overline{v}_{\varepsilon_n}(r)) \|^{2}_{_{\mathcal{L}_{2}(L^2(\mathcal S,{\mathbb R}^d) , H)}} \, \mathrm{d}r \right|^{2} \, \right] \, 
&\leq C \, \varepsilon_{n}^{2} \; \left(  \, \overline{\mathbb{E}}\left[ \, \int_{0}^{T} \|\overline{v}_{\varepsilon_n}(r)\|^{2}_{_V} \,\mathrm{d}r  \right] \, + \, 1 \right) \\
&\leq C \, \varepsilon_{n}^{2} \, .
\end{align*}
By thinning the sequence, we have also the convergence in probability 
\[
G_{\varepsilon_n}(\overline{v}_{\varepsilon_n}) \underset{n\to+\infty}{\longrightarrow} 0 \quad \text{in} \quad  L^{2}([0,T] \, ; \, \mathcal{L}_{2}(L^2(\mathcal S,{\mathbb R}^d) , V') )\,. 
\]
As in the section \ref{PassageLimite}, we apply lemma \ref{ConvSto} and the Vitali convergence theorem to obtain
\[
K_{4}^{\varepsilon_n} \underset{n\to+\infty}{\longrightarrow} 0   \qquad  \overline{\mathbb{P}} \; \text{a.s} \, .
\]
We take hence the limits when $[n\to\infty]$ in {\eqref{eq-Kin}}. By density of $\D(A^{3/2})$ in $V$, we have for almost $t\in [0,T]$ {and for almost $\omega \in \overline{\Omega}$}
\[ 
\overline{v}(t) - v_0 \, + \, \displaystyle \int_{0}^{t} \, A\,\overline{v}(r) \, \mathrm{d}r \, + \, \displaystyle \int_{0}^{t} B\, \overline{v}(r)\, \mathrm{d}r \; = \; 0
\]
in $V$ with $\overline{v} \in {L^{2}(\overline{\Omega}, L^2([0, T ] , V ) )\cap L^{2}(\overline{\Omega}, L^{\infty}([0, T ] , H))}$, $\overline{v}$ is hence 
$\overline{\mathbb{P}} $ a weak solution of the deterministic Navier-Stokes equation \eqref{NVSdeterEQabstrac}.

\newpage
\appendix
\section{Mathematical tools}
In this part, we recall some results used to prove  theorem \ref{Theorem}. 

\subsection{Compact embedding result}
The following results are proved in \cite{Flandoli-Gatarek-95} and are variations of the compactness theorems of \cite{Lions-69}, Ch. I,
Sect. 5, and \cite{Temam-83}, Sect. 13.3. 

\begin{lemma}\label{Compact1}
Let $B_{0} \subset B \subset B_{1}$ be Banach spaces with compact embedding of $B_0$ in $B$. Let $p\in (1 , \infty)$ and $\alpha \in (0,1)$ be given. Then the embedding 
\[
L^{p}([0,T] ; B_0) \cap W^{\alpha , p}([0,T ] ; B_{1}) \hookrightarrow L^{p}([0,T] ; B) \quad \text{is compact .}
\]
\end{lemma}

\begin{lemma}\label{Compact2}
Let $B_0 \subset B$ two Banach spaces with compact embedding. Let $\alpha \in (0,1)$ and $p>1$ be given such that $\alpha p >1$, then the injection 
\[
W^{\alpha , p}([0,T ] ; B_0) \hookrightarrow C^{0}([0,T] ; B) \quad \text{is compact .}
\]
\end{lemma}

\subsubsection*{Stochastic inequalities}
Let $(\Omega,\mathcal{F}, (\mathcal{F}_t)_t , \mathbb{P})$ be a stochastic basis. Let $X$ be a separable Hilbert space. Let $U$ be a second separable Hilbert space, and let $W$ be a cylindrical Wiener process with values in $U$, defined on the stochastic basis. For any progressively measurable process $\Phi \in  L^{2}(\Omega \times [0, T];\mathcal{L}_{2}(U,X))$, we denote
by $I(\Phi)$ the Ito integral defined for $t\in [0,T]$ as
\[
I(\Phi)(t) = \int_{0}^{t} \Phi(s) \, \mathrm{d} W_s . 
\]
\begin{lemma}\label{Lemme21}
Let $p\geq 2$ and $\alpha < 1/2$ be given. Then, for any progressively measurable process $\Phi \in L^{p}(\Omega \times [0, T];\mathcal{L}_{2}(U,X))$, we have
\[
I(\Phi) \in L^{p}(\Omega ; W^{\alpha ,p}([0,T] \, ; \, X))
\]
and there exists a constant $C(p,\alpha) > 0$ independent of $\Phi$ such that
\[
\mathbb{E} \left[ \, \| I(\Phi) \, \|^{p}_{W^{\alpha,p}([0,T] \,  ; \, X)}  \right]  \; \leq \; C(p, \alpha) \; \mathbb{E} \left[ \, \int_{0}^{T} \| \Phi  \|^{p}_{\mathcal{L}_{2}(U,X)} \, \mathrm{d}t \; \right] .
\]
\end{lemma}
This result is proved in \cite{Flandoli-Gatarek-95}. Most notably for the analysis here, the Burkholder-Davis-Gundy inequality holds and is recalled in the next proposition.
\begin{proposition}[Burkholder-Davis-Gundy inequality]\label{BDG} \ \\
For all integer $r\geq 1$, we have 
\[
\mathbb{E} \left[ \, \sup_{t\in[0,T]} \left| \int_{0}^{t}  \Phi(s)  \, \mathrm{d}W_s \, \right|^{r} \; \right] \; \leq \; C \, \mathbb{E} \left[ \,\left( \int_{0}^{T}  \| \, \Phi(s)\, \|^{2}_{\mathcal{L}_{2}(U,X)} \, \mathrm{d}s \;  \right)^{r/2}\right]
\]
with a constant $C>0$ depending only on $r$. \\
\end{proposition}
In section \ref{PassageLimite}, the following lemma of \cite{Debussche-et-al-2011}  is used to facilitate the passage to the limit in the Galerkin scheme. 
\begin{lemma}\label{ConvSto}
Consider a sequence of stochastic bases $S_n = (\Omega,\mathcal{F},(\mathcal{F}_{t}^{n})_{t\in[0,T]},\mathbb{P} ,W^n)$ with $W^n$ a cylindrical Wiener process (over $U$) with respect to $\mathcal{F}_{t}^{n}$. Assume that $(\psi^n)_{n\geq 1}$ is a collection of $X$-valued $\mathcal{F}_{t}^{n}$ predictable processes such that $\psi^{n} \in L^{2}\left([0,T],\mathcal{L}_2(U,X)\right)$ a.s. Finally consider $S = (\Omega,\mathcal{F}, \mathbb{P} , W)$ with $W$ a cylindrical Wiener process (over $U$) and $\psi \in L^{2}\left([0,T],\mathcal{L}_2(U,X)\right)$, which is $F_t$ predictable.
If 
\begin{align*}
W^{n}  \underset{n\to+\infty}{\longrightarrow} W \qquad \text{in probability in } \; C^{0}\left( [0,T] ; U_0 \,\right)\\
\psi^{n} \underset{n\to+\infty}{\longrightarrow} \psi \qquad \text{in probability in } \; L^{2}\left([0,T] , \mathcal{L}_{2}(U , X)\,\right) 
\end{align*} 
then 
\[
\displaystyle \int_{0}^{t} \psi^n \, \mathrm{d}W^{n}_s \underset{n\to+\infty}{\longrightarrow} \int_{0}^{t} \psi \, \mathrm{d}W_s \qquad \text{in probability in} \; L^{2}\left( [0,T] ; X \right)
\]
\end{lemma}
\subsection{Skohorod and Prokhorov theorems}
We recall here two classical theorems whose proofs can be found in \cite{dapratostochastic1992}. 
\begin{theorem}[Prokhorov]
A set $\Lambda$ of probability measures on $(E,\mathcal{B})$ is relatively compact if and only if it is tight.
\end{theorem}
\begin{theorem}[Skorohod]
\label{ThmSkorohod}
Assuming there is a sequence $\{\mu_n\}_{n\geq1}$ converging weakly to a measure $\mu$, then there exists a probability space $(\overline{\Omega},\overline{\mathcal F}, \overline{\mathbb{P}})$ and a sequence of  X-valued random variables $(\overline{Y}_n)_{n\geq0}$ (relative to this space) such that $\overline{Y}_n$ converges almost surely to the random variable $\overline{Y}$ and such that the laws of $\overline{Y}_n$ and $\overline{Y}$ are $\mu_n$ and $\mu$, respectively, i.e. $\mu_n(E) = \mathbb{P}(\overline{Y}_n \in E)$, $\mu(E) = \mathbb{P}(\overline{Y} \in E)$, for all $E \in \mathcal{B}(X)$. 
\end{theorem}

\subsection{Deterministic Navier Stokes equation }\label{NVSdeterministe}
Let $v_0\in H$. The standard deterministic Navier-Stokes equations for incompressible fluids reads
\begin{equation}\label{NVSdeterministequation}
\left\{
    \begin{array}{ll}
        d_t v(t) \; + \; \dfrac{1}{\, R_e\, } \Delta v(t)  \, \mathrm{d}t \; + \; \bigl(v(t) \bcdot \nabla \bigr) v(t)\, \mathrm{d}t \;  = \;  -\dfrac{1}{\rho} \; \nabla p \dif t, &  \\
        \nabla \bcdot v = 0, &  \\
        v(0)=v_0. & 
    \end{array}
\right.
\end{equation}
This equation can be written as the following pressure-free abstract problem using the Leray projection and the operator $A$ and $B$ defined  section $\ref{Operators}$ : 
\begin{equation}\label{EquationAbstactDeterministe}
\left\{
    \begin{array}{ll}
        d_t v(t) \; + \; Av(t)  \, \mathrm{d}t \; + B v(t)\, \mathrm{d}t \; = \; 0, &  \\
        v(0)=v_0. & 
    \end{array}
\right.
\end{equation}
We say that $v$ is a weak solution of \eqref{EquationAbstactDeterministe} if $v\in L^{2}([0,T] , V) \cap {L^\infty}([0,T] , H)$ and verifies for all $z\in V$ and $t\in [0,T]$ the following equality
\begin{equation*}
    \bigl( v(t) \, , \, z\bigr)_{\!{_H}} - \bigl( v_0 \, , \, z\bigr)_{\!{_H}} + \int_{0}^{t} \bigl( A\,v(r)\,, \, z\bigr)_{\!{_H}} \mathrm{d}r \; + \; \int_{0}^{t} \left(  B\,v(r)\, , \, z   \right)_{\!{_H}} \mathrm{d}r
\; = 0 \; .
\end{equation*}
For $d=2$ or $3$, system \eqref{EquationAbstactDeterministe} admits a weak solution. This solution is unique and {belongs to $C([0,T],H)$ } when $d=2$. 
\section*{Funding and/or Conflicts of interests/Competing interests} The author declares that there is no conflict of interest and no competing interest. The authors acknowledge the support of the ERC EU project 856408-STUOD.
\bibliographystyle{plain}
\bibliography{glob} 
\end{document}